\documentstyle[12pt]{article}
\newtheorem{Thm}{\indent Theorem}[section]
\newtheorem{Prop}[Thm]{\indent Proposition}
\newtheorem{Lem}[Thm]{\indent Lemma}
\newtheorem{Cor}[Thm]{\indent Corollary}

\newcommand{\Proof}{{\sl Proof.}\quad}
\newcommand{\QED}{{\unskip\nobreak\hfil\penalty50\quad\null\nobreak\hfil
{\sl q.e.d.}\parfillskip0pt\finalhyphendemerits0\par\medskip}}

\renewcommand{\H}{\mathop{\rm H}\nolimits}

\newcommand{\Hom}{\mathop{\rm Hom}\nolimits}

\newcommand{\IC}{\mathop{\rm IC}\nolimits}
\newcommand{\GM}{\mathop{\rm GM}\nolimits}
\newcommand{\CGM}{\mathop{\rm CGM}\nolimits}
\newcommand{\GEM}{\mathop{\rm GEM}\nolimits}
\newcommand{\UGEM}{\mathop{\rm UGEM}\nolimits}
\newcommand{\CGEM}{\mathop{\rm CGEM}\nolimits}
\newcommand{\Coh}{\mathop{\rm Coh}\nolimits}
\newcommand{\CCohDiff}{\mathop{\rm CCohDiff}\nolimits}

\newcommand{\Det}{\mathop{\rm Det}\nolimits}
\newcommand{\ic}{\mathop{\rm ic}\nolimits}

\newcommand{\grt}{\mathop{\rm gt}\nolimits}
\newcommand{\relint}{\mathop{\rm rel.\,int}\nolimits}

\newcommand{\setmin}{{\setminus}}
\newcommand{\inte}{\mathop{\rm int}\nolimits}
\newcommand{\Ker}{\mathop{\rm Ker}\nolimits}
\newcommand{\Coker}{\mathop{\rm Coker}\nolimits}
\renewcommand{\Im}{\mathop{\rm Im}\nolimits}
\newcommand{\Spec}{\mathop{\rm Spec}\nolimits}

\newcommand{\rank}{\mathop{\rm rank}\nolimits}

\newcommand{\scl}{\mathrel{;}}
\newcommand{\cl}{\mathrel{:}}

\newcommand{\Zero}{{\bf 0}}
\newcommand{\Gm}{{\bf G}_{\rm m}}
\newcommand{\calC}{{\cal C}}

\newcommand{\calF}{{\cal F}}
\newcommand{\calH}{{\cal H}}

\newcommand{\calO}{{\cal O}}
\newcommand{\calP}{{\cal P}}

\newcommand{\R}{{\bf R}}
\newcommand{\Z}{{\bf Z}}
\newcommand{\C}{{\bf C}}
\renewcommand{\d}{{\rm D}}
\newcommand{\D}{{\bf D}}
\newcommand{\Q}{{\bf Q}}
\renewcommand{\i}{{\rm i}}
\renewcommand{\j}{{\rm j}}
\newcommand{\e}{{\bf e}}
\newcommand{\m}{{\bf m}}
\newcommand{\p}{{\bf p}}
\renewcommand{\t}{{\bf t}}
\newcommand{\h}{{\rm h}}
\newcommand{\lan}{\langle}
\newcommand{\ran}{\rangle}
\newcommand{\lra}{\longrightarrow}
\newcommand{\ra}{\rightarrow}

\title
{Torus embeddings and\\ algebraic intersection complexes}

\author
{Masa-Nori ISHIDA\thanks{Supported in part by a Grant under
The Monbusho International Scientific Research Program: 04044081}
\\
{\normalsize Mathematical Institute, Faculty of Science}\\
{\normalsize Tohoku University, Sendai 980, Japan}}

\begin{document}

\maketitle%


\section*{Introduction}




In \cite{GM2}, Goreskey and MacPherson defined and constructed
intersection complexes for topological pseudomanifolds.
The complexes are defined in the derived category of sheaves
of modules over a constant ring sheaf.
Since analytic spaces are of this category, any algebraic variety
defined over $\C$ has an intersection complex for each perversity.

The purpose of this paper is to give an algebraic description of the
intersection complex of a toric variety.
Namely, we describe it as a finite complex of coherent sheaves
whose coboudary map is a differential operator of order one.

Let $Z_\h$ be the complete toric variety associated to a complete
fan $\Delta$.
For each $\sigma\in\Delta$, let $X(\sigma)_\h$ be the associated
closed subvariety of $Z_\h$.
For each perversity $\p$, we construct a bicomplex
$(\ic_\p(Z_\h)^{\bullet,\bullet}, d_1, d_2)$ with the following
properties.

(1) $\ic_\p(Z_\h)^{i,j} = \{0\}$ for $(i,j)\not\in[0,r]\times[-r,0]$.

(2) Each $\ic_\p(Z_\h)^{i,j}$ is a direct sum for $\sigma\in\Delta$
of free $\calO_{X(\sigma)_\h}$-modules of finite rank.

(3) $d_1$ is an $\calO_{Z_\h}$-homomorphism and $d_2$ is a differential
operator of order one.

(4) The associated single complex $\ic_\p(Z_\h)^\bullet$ is
quasi-isomorphic to the $r$-times dimension shifts to the right of
the intersection complex defined in \cite{GM2}.
In other words, our complex belongs to
``Beilinson-Bernstein-Deligne-Gabber scheme'' (cf. \cite[2.3,(d)]{GM2}
and \cite[2.1]{BBD}).

In \S1, we introduce an abelian category $\GM(A(\sigma))$ of
finitely generated graded $A(\sigma)$-modules, where $A(\sigma)$
is the exterior algebra of the $\Q$-vector space $N(\sigma)_\Q$
defined by a cone $\sigma$.

In \S2, we define an additive category $\GEM(\Delta)$ for
a finite fan $\Delta$.
Each object $L$ of this category is a collection
$(L(\sigma)\scl\sigma\in\Delta)$ of $L(\sigma)\in\GM(A(\sigma))$.
We define a dualizing functor $\D$ on the category
$\CGEM(\Delta)$ of finite complexes in $\GEM(\Delta)$.

A perversity $\p$ on $\Delta$ is defined to be a
$\Z$-valued map on $\Delta\setminus\{\Zero\}$.
The intersection complex $\ic_\p(\Delta)^\bullet$ as
an object of $\CGEM(\Delta)$ is defined and constructed
in \S2.

In \S3, we work on the toric variety $Z(\Delta)$ associated
to the fan $\Delta$.
For each $L^\bullet\in\CGEM(\Delta)$, we define a finite
bicomplex $\Lambda_{Z(\Delta)}(L)^{\bullet,\bullet}$ of
coherent $\calO_{Z(\Delta)}$-modules whose second
coboundary map is a differential operator of order one.

We consider the normal analytic space $Z_\h := Z(\Delta)_\h$
associated to the toric variety $Z(\Delta)$ in \S4.
The bicomplex $\ic_\p(Z_\h)^{\bullet,\bullet}$
stated above is defined to be the bicomplex
$\Lambda_{Z(\Delta)_\h}(\ic_\p(\Delta))^{\bullet,\bullet}$
on this analytic space.

When $\Delta$ is a complete fan, we show that the
intersection cohomologies are described in termes of the
complex $\ic_\p(\Delta)^\bullet$ in $\CGEM(\Delta)$.

The middle perversity $\m$ is defined by $\m(\sigma) := 0$
for all $\sigma\in\Delta\setminus\{\Zero\}$.
In \cite{Ishida3}, we will discuss on $\ic_\m(\Delta)^\bullet$
and prove the decomposition theorem for a barycentric
subdivision of the fan.


\section*{Notation}




We denote by $\Z$ the ring of rational integers and by $\Q$, $\R$
and $\C$ the fields of rational numbers, real numbers and
complex numbers, respectively.

For a free $\Z$-module $F$ of finite rank, we denote
$F_\Q := F\otimes\Q$ and $F_\R := F\otimes\R$.

We denote a complex $E^\bullet$ of modules or of sheaves of
modules simply by $E$ when the substitution of the dot by
an integer is not suitable.
In particular, we prefer to write $G(E)^\bullet$ rather than
$G(E^\bullet)$, if $G$ is a functor between categories of complexes.
However, we left the dot in representing the cohomologies
$\H^i(E^\bullet)$.

For a complex $(E^\bullet, d_E)$ and a integer $n$, the complex
$(E[n]^\bullet, d_{E[n]})$ is defined by $E[n]^i := E^{i+n}$ and
$d_{E[n]}^i := (-1)^nd_E^{i+n}$ for $i\in\Z$.

By a bicomplex $(E^{\bullet,\bullet}, d_1, d_2)$, we mean a
naive double complex \cite[0.4]{Deligne}, i.e., a double complex
satisfying $d_1\cdot d_2 = d_2\cdot d_1$.
The associated single complex $(E^\bullet, d)$ is defined by
$E^k :=\bigoplus_{i+j=k}E^{i,j}$ for every $k\in\Z$ and
$d(x) := d_1^{i,j}(x) + (-1)^id_2^{i,j}(x)$ for
$(i, j)\in\Z\times\Z$ and $x\in E^{i,j}$.
In general, we follow \cite{Deligne} for the sign
configuration of complexes.

By a $d$-complex and a $\partial$-complex, we mean complexes in
an additive category whose coboundary
maps are denoted by $d$ and $\partial$, respectively.
Some impotant bicomplexes in this paper has $d_1 = d$ and
$d_2 =\partial$.


\section{The exterior algebras and modules}



\setcounter{equation}{0}

Let $r$ be a non-negative integer, and let $N$ and $M$ be
mutually dual free $\Z$-modules of rank $r$ with the pairing
$\lan\;,\;\ran\cl M\times N\lra\Z$.
This pairing is extended $\R$-bilinearly to
$\lan\;,\;\ran\cl M_\R\times N_\R\lra\R$.

By a {\em cone} in $N_\R$, we mean a strongly convex rational
polyhedral cone \cite[Chap.1,1.1]{Oda1},
i.e., a cone $\sigma$ is equal to
$\R_0n_1 +\cdots +\R_0n_s$ for a finite subset
$\{n_1,\cdots, n_s\}$ of $N$ and satisfies
$\sigma\cap(-\sigma) =\{0\}$, where
$\R_0 :=\{c\in\R\scl c\geq 0\}$.
The cone $\{0\}$ is denoted by $\Zero$.
For each cone $\sigma$ we set $r_\sigma :=\dim\sigma$.

For a cone $\sigma$ in $N_\R$, we define
$N(\sigma) := N\cap(\sigma + (-\sigma))\simeq\Z^{r_\sigma}$ and
$N[\sigma] := N/N(\sigma)\simeq\Z^{r-r_\sigma}$.
Hence, $N(\sigma)_\R$ is the real subspace $\sigma + (-\sigma)$
of $N_\R$.
On the other hand, we set
$\sigma^\bot :=\{x\in M_\R\scl\lan x, a\ran = 0,\forall a\in\sigma\}$,
$M[\sigma] := M\cap\sigma^\bot$ and $M(\sigma) := M/M[\sigma]$.
Then $N(\sigma)_\R\subset N_\R$ and $M[\sigma]_\R\subset M_\R$
are orthogonal complement of each other with respect to the pairing.
These notations are defined so that $N(\sigma)$ and $M(\sigma)$
as well as $N[\sigma]$ and $M[\sigma]$ are mutually dual,
respectively.

In this article, we often treat finite dimensional
graded $\Q$-vector spaces.
Let $V$ be such a $\Q$-vector space.
Then, we denote by $V_j$ the vector subspace consisting of zero
and the homogeneous elements of degree $j$ for each $j\in\Z$.
For two finite dimensional $\Q$-vector spaces $V$ and
$W$, we define the grading of $V\otimes_\Q W$ by
\begin{equation}
(V\otimes_\Q W)_k :=\bigoplus_{i+j=k}V_i\otimes_\Q W_j
\end{equation}
for each $k\in\Z$.
We identify $V\otimes_\Q W$ with $W\otimes_\Q V$ by the
identifications
\begin{equation}
\label{left and right op}
a_i\otimes b_j = (-1)^{ij}b_j\otimes a_i
\end{equation}
for $a_i\in V_i, b_j\in W_j$ for all $i, j\in\Z$.

We denote by $A(M_\Q)$ the exterior algebra
$\bigwedge^\bullet M_\Q$ over the rational number field $\Q$.
Then $A(M_\Q)$ is a graded $\Q$-algebra with
$A(M_\Q)_i :=\bigwedge^i M_\Q$ for each $i\in\Z$.

The algebra $A(N_\Q) :=\bigwedge^\bullet N_\Q$ is more important
in our theory.
The grading of $A(N_\Q)$ is defined negatively, i.e.,
$A(N_\Q)_i :=\bigwedge^{-i} N_\Q$ for each $i$.
For a cone $\sigma$ in $N_\R$,
$A(N(\sigma)_\Q) :=\bigwedge^\bullet N(\sigma)_\Q$ is a graded
subalgebra of $A(N_\Q)$.

In order to simplify the notation, we set $A := A(N_\Q)$ and
$A^* := A(M_\Q)$.
For a cone $\sigma$ in $N_\R$, we set
$A(\sigma) := A(N(\sigma)_\Q)\subset A$
and $A^*[\sigma] := A(M[\sigma]_\Q)\subset A^*$.

Let $C =\bigoplus_{p\in\Z}C_p$ be a graded $\Q$-subalgebra
of $A$ or $A^*$.
For a graded left $C$-module $V =\bigoplus_{q\in\Z}V_q$, we define
the {\em associated} graded right $C$-module structure on $V$ by
\begin{equation}
\label{leftright}
xa := \sum_{q\in\Z}(-1)^{pq}ax_q
\end{equation}
for $a\in C_p$ and $x =\sum_{q\in\Z}x_q\in V$.
Conversely, if $V$ is a graded right $C$-module, then the asssociated
graded left $C$-module structure is defined similarly.
In both cases, we see easily that $a(xb) = (ax)b$ for $a, b\in C$
and $x\in V$.
Hence $V$ is a two-sided $C$-module and we call it simply a
$C$-module.
The following lemma is checked easily.

%
%

\begin{Lem}
\label{Lem 1.1}
Let $V, V'$ be graded $C$-modules for $C$ as above.
If $W$ is a homogeneous left or right $C$-submodule of $V$,
then it is a two-sided $C$-submodule of $V$.
If $f\cl V\ra V'$ is a homogeneous homomorphism as left or
right $C$-modules of degree zero, then it is a homomorphism
of two-sided $C$-modules.
\end{Lem}

Let $\GM(C)$ be the abelian category of finitely generated
graded $C$-modules where the morphisms are defined to be
homogeneous homomorphisms of degree zero.
By definition, every object in $\GM(C)$ is
finite-dimensional as a $\Q$-vector space.

Since $M_\Q$ is the dual $\Q$-vector space of $N_\Q$, each
homogeneous element $a\in A^*_p$ for any $p\in\Z$ induces a
homogeneous $\Q$-linear map $i(a)\cl A\ra A$ (cf.\ \cite[5.14]{Greub})
which is usually called the right interior product.
Here note that the degree of the interior product $i(a)$ is $p$,
since the indices of $A$ are given negatively.
It is known that this operation induces a graded right
$A^*$-module structure on $A$ (cf.\ \cite[(5.50)]{Greub}).
With the associated left $A^*$-module structure, we regard
$A$ as a two-sided $A^*$-module.

%
%

\begin{Lem}
\label{Lem 1.2}
Let $\sigma$ be a cone in $N_\R$.
Then the left operations of $A(\sigma)$ and
$A^*[\sigma]$ on $A$ commute with each other.
This is also true for the right operations.
\end{Lem}

\Proof
For $a\in A^*[\sigma]_p$, $b\in A(\sigma)_q$ and
$x\in A_s$, we prove the equalities
\begin{eqnarray}
b(ax) & = & a(bx)\;, \label{Lem 1.2 1} \\
(xa)b & = & (xb)a\;, \label{Lem 1.2 2} \\
(bx)a & = & (-1)^{pq}b(xa)\;. \label{Lem 1.2 3}
\end{eqnarray}

Since $N(\sigma)_\Q$ is the orthogonal complement of
$M[\sigma]_\Q$ and $(bx)a = i(a)(b\wedge x)$, the equality
(\ref{Lem 1.2 3}) is equal to \cite[(5.58)]{Greub}.
By this equality, we have
$b(ax) = (-1)^{ps}b(xa) = (-1)^{ps+pq}(bx)a = a(bx)$ and
$(xa)b = (-1)^{pq+qs}b(xa) = (-1)^{qs}(bx)a =(xb)a$.
These are the equalities (\ref{Lem 1.2 1}) and (\ref{Lem 1.2 2}).
\QED

Let $\sigma$ be a cone in $N_\R$ and let $V$ be a graded
$A(\sigma)$-module.
Then the above lemma implies that $V_A := V\otimes_{A(\sigma)}A$
has a structure of $A^*[\sigma]$-module such that
$a(u\otimes x) = u\otimes ax$ for $u\in V$, $x\in A$ and
$a\in A^*[\sigma]$.

%
%

\begin{Lem}
\label{Lem 1.3}
Let $V$ be a finitely generated graded $A(\sigma)$-module
with a homogeneous $\Q$-basis $\{u_1,\cdots, u_s\}$ and let
$\{x_1,\cdots, x_r\}$ be a $\Q$-basis of $N_\Q$ such that
$N(\sigma)_\Q$ is generated by $\{x_1,\cdots, x_k\}$  for
$k :=r_\sigma$.
Then $V_A$ is a free $A^*[\sigma]$-module with the
basis $\{u'_1,\cdots, u'_s\}$, where
$u'_i :=u_i\otimes(x_{k+1}\wedge\cdots\wedge x_r)$ for each $i$.
\end{Lem}

\Proof
Let $E$ be the subspace $\Q x_{k+1} +\cdots +\Q x_r$ of $N_\Q$.
Since the operation of $A^*[\sigma]$ on $A$ is defined
by the interior products, we have
$A(E) = A^*[\sigma](x_{k+1}\wedge\cdots\wedge x_r)$ in $A$.
Hence $A = A(\sigma)\otimes_\Q A(E)$ is equal to
$(A(\sigma)\otimes_\Q A^*[\sigma])(x_{k+1}\wedge\cdots\wedge x_r)$.
Hence
\begin{equation}
V_A = V\otimes_{A(\sigma)}A\simeq
V\otimes_\Q A^*[\sigma](x_{k+1}\wedge\cdots\wedge x_r)
\end{equation}
as $A^*[\sigma]$-modules and we get the lemma.
\QED

Let $\sigma$ be a cone of $N_\R$.
We set
\begin{equation}
\det(\sigma) := \bigwedge^{r_\sigma}N(\sigma)\simeq\Z
\end{equation}
and $\det(\sigma)_\Q :=\det(\sigma)\otimes\Q$.
We denote by $\Det(\sigma)_\Q$ the graded $\Q$-vector space
defined by $(\Det(\sigma)_\Q)_{-r_\sigma} :=\det(\sigma)_\Q$
and $(\Det(\sigma)_\Q)_j :=\{0\}$ for $j\not = -r_\sigma$.

For a finitely generated graded $A(\sigma)$-module $V$,
we define a graded $A(\sigma)$-module $\d_\sigma(V)$ as follows.

We define graded $\Q$-vector spaces $\d_\sigma^{\rm left}(V)$ and
$\d_\sigma^{\rm right}(V)$ by
\begin{equation}
\d_\sigma^{\rm left}(V) =\d_\sigma^{\rm right}(V) :=
\Hom_\Q(V,\Det(\sigma)_\Q) =
\bigoplus_{i\in\Z}\Hom_\Q(V_{-r_\sigma-i},\det(\sigma)_\Q)\;.
\end{equation}
For $x\in V$ and $y\in\d_\sigma^{\rm left}(V)$
(resp.\ $z\in\d_\sigma^{\rm right}(V)$) we denote the operation by
$(y, x)$ (resp.\ (x, z)).
The right $A(\sigma)$-module structure of $\d_\sigma^{\rm left}(V)$
(resp.\ the left $A(\sigma)$-module structure of
$\d_\sigma^{\rm right}(V)$) is defined by
\begin{equation}
(ya, x) := (y, ax)\;\;\;(\hbox{ resp.\ }(x, az) := (xa, z))
\end{equation}
for $a\in A(\sigma)$.
There exists a unique homogeneous isomorphism
$\varphi\cl\d_\sigma^{\rm left}(V)\ra\d_\sigma^{\rm right}(V)$
degree zero such that $(y, x) = (-1)^{pq}(x, \varphi(y))$
for homogeneous elements $x\in V_p$ and $y\in\d_\sigma^{\rm left}(V)_q$
for all $p, q\in\Z$.
We define $\d_\sigma(V)$ to be the identification of
$\d_\sigma^{\rm left}(V)$ and $\d_\sigma^{\rm right}(V)$ by the
isomorphism $\varphi$.
Namely, we have $(y, x) = (-1)^{pq}(x, y)$ for $x\in V_p$ and
$y\in\d_\sigma(V)_q$.
Here note that $(y, x) = 0$ if $p + q\not= -r_\sigma$.
It is easy to see that the induced left and right
$A(\sigma)$-module structures on $\d_\sigma(V)$ have the
compatibility (\ref{leftright}).
By definition, we have
\begin{equation}
\label{eq of d-dual}
\dim_\Q\d_\sigma(V)_j =\dim_\Q V_{-r_\sigma-j}
\end{equation}
for every $j\in\Z$.

For $V\in\GM(A(\sigma))$, we define an $A(\sigma)$-homomorphism
$\iota\cl V\ra\d_\sigma(\d_\sigma(V))$ by
$(y,\iota(x)) := (y, x)$ for $x\in V$ and $y\in\d_\sigma(V)$.
It is easy to see that the symmetric equality
$(\iota(x), y) = (x, y)$ holds.
Since $V$ is a finite dimensional $\Q$-vector space, the
pairings are perfect and $\iota$ is an isomorphism.
We call $\iota$ the {\em canoniacal isomorphism}.

For a homomorphism $f\cl V\ra W$ in $\GM(A(\sigma))$,
we define $\d_\sigma(f)$ to be the natural induced homomorphism
$\d_\sigma(W)\ra \d_\sigma(V)$ of $A(\sigma)$-modules,
i.e., the equality $(y, f(x)) = (\d_\sigma(f)(y), x)$ for
$x\in V$ and $y\in\d_\sigma(W)$.
It is clear by definition that the correspondence $V\mapsto\d_\sigma(V)$
is a contravariant exact functor from $\GM(A(\sigma))$ to itself.

If $\pi$ is of dimension $r$, then $A(\pi) = A$ and
$\d_\pi$ is a functor from $\GM(A)$ to itself.
We denote this functor by $\d_N$ which does not depend
on the choice of $\pi$.

Let $\sigma$ and $\rho$ be cones in $N_\R$ with $\sigma\prec\rho$.

For $V$ in $\GM(A(\sigma))$, we denote by $V_{A(\rho)}$ the
graded $A(\rho)$-module
$V\otimes_{A(\sigma)}A(\rho) = A(\rho)\otimes_{A(\sigma)}V$,
where we identify $x\otimes a$ with $(-1)^{pq}a\otimes x$ for
$x\in V_p$ and $a\in A(\rho)_q$.
For a morphism $f\cl V\ra V'$ in $\GM(A(\sigma))$, we denote
$f_{A(\rho)} := f\otimes1_{A(\rho)}\cl V_{A(\rho)}\ra V'_{A(\rho)}$.
The correspondence $V\mapsto V_{A(\rho)}$ is a covariant functor from
$\GM(A(\sigma))$ to $\GM(A(\rho))$.
Let $H$ be a linear subspace of $N(\rho)_\Q$ such that
$N(\rho)_\Q = N(\sigma)_\Q\oplus H$.
Then $A(\rho) = A(\sigma)\otimes_\Q A(H)$ and
$V_{A(\rho)} = V\otimes_\Q A(H)$ for any $V$ in $\GM(A(\sigma))$.
This implies that the functor is exact.
Similarly, the functor from $\GM(A(\sigma))$ to $\GM(A)$ defined
by $V\mapsto V_A$ is exact.

For $V\in\GM(A(\sigma))$, we define an $A(\sigma)$-homomorphism
\begin{equation}
\label{d_ sigma rho}
\varphi\cl\d_\sigma(V)\ra\d_\rho(V_{A(\rho)})
\end{equation}
by $(\varphi(y), xa) := \phi_\rho((y, x)a)$ for $x\in V$,
$y\in\d_\sigma(V)$ and $a\in A(\rho)$, where $\phi_\rho$ is the
homogeneous projection $A(\rho)\ra A(\rho)_{-r_\rho} =\Det(\rho)_\Q$.

%
%

\begin{Lem}
\label{Lem 1.4}
Let $\sigma,\rho$ be cones in $N_\R$ with $\sigma\prec\rho$.
Then the homomorphism (\ref{d_ sigma rho}) induces $A(\rho)$-isomorphism
$\d_\rho(V_{A(\rho)})\simeq\d_\sigma(V)_{A(\rho)}$ for every $V$ in
$\GM(A(\sigma))$.
\end{Lem}

\Proof
Consider the case $V = A(\sigma)$.
Let $\phi_\sigma\cl A(\sigma)\ra\Det(\sigma)_\Q = A(\sigma)_{-r_\sigma}$
be the homogeneous projection.
For $u\in A(\sigma)$, the corresponding  element $\phi_\sigma u$ of
$\d_\sigma(A(\sigma))$ is given by
$(\phi_\sigma u, x) :=\phi_\sigma(u\wedge x)$ for $x\in A(\sigma)$.
Hence $\d_\sigma(A(\sigma))$ is a free $A(\sigma)$-module
generated by $\phi_\sigma$.
Similarly, $\d_\rho(A(\rho))$ is equal to $\phi_\rho A(\rho)$.

By the definiton of $\varphi$ in (\ref{d_ sigma rho}), we have
$\varphi(\phi_\sigma) = \phi_\rho$.
Hence $\varphi$ induces an isomorphism
$\d_\sigma(A(\sigma))_{A(\rho)}\simeq\d_\rho(A(\rho))$.

For general $V$, we take an exact sequence
\begin{equation}
\bigoplus_{i=1}^m A(\sigma)x_i\mathop{\lra}\limits^{f}
\bigoplus_{j=1}^n A(\sigma)y_j\lra V\lra0
\end{equation}
of graded left $A(\sigma)$-modules, where $\{x_1,\cdots, x_m\}$
and $\{y_1,\cdots, y_n\}$ are homogeneous bases.
Then by the exactness of the functors, we get exact sequences
\begin{equation}
0\lra\d_\sigma(V)_{A(\rho)}\lra\bigoplus_{j=1}^n y_j^*A(\rho)
\mathop{\lra}\limits^{{}^{\rm t}f}\bigoplus_{i=1}^mx_i^*A(\rho)
\end{equation}
and
\begin{equation}
0\lra\d_\rho(V_{A(\rho)})\lra\bigoplus_{j=1}^n y_j^*A(\rho)
\mathop{\lra}\limits^{{}^{\rm t}f}\bigoplus_{i=1}^mx_i^*A(\rho)
\end{equation}
of graded right $A(\rho)$-modules.
Hence we get the required isomorphism.
\QED

Let $V^\bullet$ be a finite $d$-complex of graded $A(\sigma)$-modules
in $\GM(A(\sigma))$ with $d = (d_V^i\cl i\in\Z)$.
We define the complex $\d_\sigma(V)^\bullet$ by
$\d_\sigma(V)^i := \d_\sigma(V^{-i})$ for $i\in\Z$.
The coboudary map $d = (d_{\d_\sigma(V)}^i)$ is defined by
\begin{equation}
d_{\d_\sigma(V)}^i := (-1)^{i+1}\d_\sigma(d_V^{-i-1})
\cl\d_\sigma(V)^i\lra\d_\sigma(V)^{i+1}
\end{equation}
for each $i\in\Z$ (cf. \cite[1.1.5]{Deligne}).
Note that we have the principle to put the dot sign of complexes
at the right end.

Since $\d_\sigma$ is an exact functor, the cohomology group
$\H^p(\d_\sigma(V)^\bullet)$ is isomorphic to
$\d_\sigma(\H^{-p}(V^\bullet))$ as a graded $\Q$-vector space.
By the equality (\ref{eq of d-dual}), we get the following lemma.

%
%

\begin{Lem}
\label{Lem 1.5}
Let $\sigma\subset N_\R$ be a cone and let
$V^\bullet$ be a finite $d$-complex in $\GM(A(\sigma))$.
Then
\begin{equation}
\dim_\Q\H^p(\d_\sigma(V)^\bullet)_q =
\dim_\Q\H^{-p}(V^\bullet)_{-r_\sigma-q}
\end{equation}
for any integers $p,q$.
In particular, we have
\begin{equation}
\dim_\Q\H^p(\d_N(V)^\bullet)_q =\dim_\Q\H^{-p}(V^\bullet)_{-r-q}
\end{equation}
if $V^\bullet$ is a finite $d$-complex in $\GM(A)$.
\end{Lem}

For a homomorphism $f\cl V^\bullet\ra W^\bullet$ in $\CGM(A(\sigma))$,
the homomorphism
\begin{equation}
\d_\sigma(f)\cl\d_\sigma(W)^\bullet\ra\d_\sigma(V)^\bullet
\end{equation}
is defined as the collection
\begin{equation}
\{\d_\sigma(f)^i =\d_\sigma(f^{-i})\scl i\in\Z\}\;.
\end{equation}

For a homogeneous $\Q$-subalgebra $C$ of $A$,
we denote by $\CGM(C)$ the category of finite $d$-complexes
in $\GM(C)$.
It is easy to see that $\d_\sigma$ is a contravariant exact
functor of the abelian category $\CGM(A(\sigma))$ to itself.

Let $V^\bullet$ be an object of $\CGM(C)$.
Hence each $V^i =\bigoplus_{j\in\Z}V_j^i$ is in $\GM(C)$.
For each $i, j\in\Z$, we denote by
$d^{i,j}\cl V_j^i\ra V_j^{i+1}$ the homogeneous component of
$d^i$ of degree $j$.
For each integer $k$, the {\em gradual truncation} below
$(\grt_{\leq k}V)^\bullet$ is the homogeneous subcomplex
of $V^\bullet$ defined by
\begin{equation}
(\grt_{\leq k}V)_j^i :=\left\{
\begin{array}{lll}
V_j^i    & \hbox{ if } & i+j < k \\
\Ker d^{i,j} & \hbox{ if } & i+j = k \\
\{0\}    & \hbox{ if } & i+j > k
\end{array}
\right.
\end{equation}
and the gradual truncation above
$(\grt^{\geq k}V)^\bullet$ is the homogeneous quotient
complex of $V^\bullet$ defined by
\begin{equation}
(\grt^{\geq k}V)_j^i :=\left\{
\begin{array}{lll}
\{0\}    & \hbox{ if } & i+j < k \\
\Coker d^{i-1,j} & \hbox{ if } & i+j = k \\
V_j^i    & \hbox{ if } & i+j > k\;.
\end{array}
\right.
\end{equation}
Since $C$ is graded negatively, $(\grt_{\leq k}V)^\bullet$
and $(\grt^{\geq k}V)^\bullet$ are $d$-complexes in $\GM(C)$.
Hence, these are covariant functors from $\CGM(C)$
to itself.

The variant gradual truncations $(\widetilde{\grt}_{\leq k}V)^\bullet$
and $(\widetilde{\grt}^{\geq k}V)^\bullet$ are defined by
\begin{equation}
(\widetilde{\grt}_{\leq k}V)_j^i :=\left\{
\begin{array}{lll}
V_j^i    & \hbox{ if } & i+j\leq k \\
\Im d^{i-1,j} & \hbox{ if } & i+j = k+1 \\
\{0\}    & \hbox{ if } & i+j > k+1
\end{array}
\right.
\end{equation}
and
\begin{equation}
(\widetilde{\grt}^{\geq k}V)_j^i :=\left\{
\begin{array}{lll}
\{0\}    & \hbox{ if } & i+j < k-1 \\
\Im d^{i,j} & \hbox{ if } & i+j = k-1 \\
V_j^i    & \hbox{ if } & i+j\geq k\;,
\end{array}
\right.
\end{equation}
respectively.

It is easy to see that $(\widetilde{\grt}_{\leq k}V)^\bullet$ is
quasi-isomorphic to $(\grt_{\leq k}V)^\bullet$, while
$(\widetilde{\grt}^{\geq k}V)^\bullet$ is quasi-isomorphic to
$(\grt^{\geq k}V)^\bullet$ (cf.\cite[p.93]{GM2}).

It is clear that
\begin{equation}
\label{H of gt leq}
\H^p((\grt_{\leq k}V)^\bullet)_q =\left\{\begin{array}{ll}
  \H^p(V^\bullet)_q & \hbox{ for } p + q\leq k \\
  \{0\}             & \hbox{ for } p + q > k
\end{array}\right.
\end{equation}
and
\begin{equation}
\H^p((\grt^{\geq k}V)^\bullet)_q =\left\{\begin{array}{ll}
  \{0\}             & \hbox{ for } p + q < k \\
  \H^p(V^\bullet)_q & \hbox{ for } p + q\geq k\;.
\end{array}\right.
\end{equation}

An object $V^\bullet$ in $\CGM(C)$ is said to be
{\em acyclic} if $\H^p(V^\bullet) = \{0\}$ for all integers $p$.

%
%

\begin{Lem}
\label{Lem 1.6}
Let $\sigma$ be a cone of $N_\R$.
Let $V^\bullet$ be in $\CGM(A(\sigma))$ and $k$ be an integer.
Then the $d$-complex $(\grt_{\leq k}\d_\sigma(V))^\bullet$ is acyclic
if and only if $(\grt^{\geq -r_\sigma -k}V)^\bullet$ is.
Similarly, the $d$-complex $(\grt^{\geq k}\d_\sigma(V))^\bullet$ is
acyclic if and only if $(\grt_{\leq -r_\sigma -k}V)^\bullet$ is.

If $V^\bullet$ is in $\CGM(A)$, then these assertions with
$\d_\sigma$ replaced by $\d_N$ and $r_\sigma$ replaced by
$r$ hold.
\end{Lem}

\Proof
By (\ref{H of gt leq}),
$(\grt_{\leq k}\d_\sigma(V))^\bullet$ is acyclic if and only if
$\H^p(\d_\sigma(V)^\bullet)_q =\{0\}$ for $p + q\leq k$.
By Lemma~\ref{Lem 1.5}, this is equivalent to the condition
\begin{equation}
\dim_\Q\H^p(V^\bullet)_q =
\dim_\Q\H^{-p}(\d_\sigma(V)^\bullet)_{-r_\sigma-q} = 0
\end{equation}
for $p + q\geq -r_\sigma - k$.
This condition means that $(\grt^{\geq-r_\sigma-k}V)^\bullet$ is acyclic.
The second assertion is similarly proved.

The last assertion is obtained by taking $\sigma$ with
$r_\sigma = r$.
\QED


\section{The graded exterior modules on a fan}



\setcounter{equation}{0}

Let $\Delta$ be a finite fan of $N_\R$ \cite[1.1]{Oda1}.
We introduce an additive category $\GEM(\Delta)$ which contains
$\GM(A(\sigma))$ as full subcategories for all $\sigma\in\Delta$.

Let $V$ be in $\GM(A(\sigma))$ and $W$ in $\GM(A(\rho))$.
If $\sigma\prec\rho$, then $A(\sigma)\subset A(\rho)$ and
$W$ has an induced structure of graded $A(\sigma)$-module.

A morphism $f\cl V\ra W$ in $\GEM(\Delta)$
is defined to be a homogeneous $A(\sigma)$-homomorphism
of degree zero.
If $\sigma$ is not a face of $\rho$, we allow only the zero
map as a morphism even if $A(\sigma)$ happens to be contained
in $A(\rho)$.

Consequently, the additive category $\GEM(\Delta)$ of
{\em graded exterior modules} on $\Delta$ is defined as follows.

A graded exterior module $L$ on $\Delta$ is a collection
$(L(\sigma)\scl\sigma\in\Delta)$ of objects $L(\sigma)$
in $\GM(A(\sigma))$ for $\sigma\in\Delta$.
A {\em homomorphism} $f\cl L\ra K$ of graded exterior modules
on $\Delta$ is a collection $f = (f(\sigma/\rho))$ of morphisms
\begin{equation}
f(\sigma/\rho)\cl L(\sigma)\lra K(\rho)
\end{equation}
in $\GM(A(\sigma))$ for all pairs $(\sigma,\rho)$ of cones
in $\Delta$ with $\sigma\prec\rho$.
For $f\cl L\ra K$ and $g\cl K\ra J$, the composite
$(g\cdot f)\cl L\ra J$ is defined by
\begin{equation}
(g\cdot f)(\sigma/\rho) :=
\sum_{\tau\in F[\sigma,\rho]}g(\tau/\rho)\cdot f(\sigma/\tau)
\end{equation}
for $\sigma,\rho$ with $\sigma\prec\rho$, where
$F[\sigma,\rho]$ is the set of the faces $\tau$ of $\rho$
with $\sigma\prec\tau$.

The direct sum of finite objects in $\GEM(\Delta)$ is defined
naturally.
An object $V$ of $\GM(A(\sigma))$ is also regarded as an object of
$\GEM(\Delta)$ by defining $V(\sigma) := V$ and $V(\rho) := \{0\}$
for $\rho\not= \sigma$.
In this sense, we may write
$L =\bigoplus_{\sigma\in\Delta}L(\sigma)$.

A homomorphism $f\cl L\ra K$ is said to be {\em unmixed} if
$f(\sigma/\rho) = 0$ for any $\sigma,\rho$ with $\sigma\not= \rho$.
If $f$ is unmixed, $\Ker f$, $\Coker f$ and $\Im f$ are defined
naturally as an object of $\GEM(\Delta)$.
We denote by $\UGEM(\Delta)$ the category of the objects
of $\GEM(\Delta)$ with the class of homomorphisms restricted
to unmixed ones.
It is easy to see that $\UGEM(\Delta)$ is an abelian category.

Let $f\cl L\ra K$ be a homomorphism in $\GEM(\Delta)$.
We say that $L$ is a submodule or a subobject of $K$,
if $f$ is unmixed and $f(\sigma/\sigma)\cl L(\sigma)\ra K(\sigma)$
is an inclusion map for every $\sigma\in\Delta$.
If $L$ is a submodule of $K$, then we define an object $K/L$ in
$\GEM(\Delta)$ by $(K/L)(\sigma) := K(\sigma)/L(\sigma)$ for
$\sigma\in\Delta$.
Namely, we have a short exact sequence
\begin{equation}
0\lra L\lra K\lra K/L\lra 0
\end{equation}
in $\UGEM(\Delta)$.

We denote $\hat\Delta :=\Delta\cup\{\alpha\}$ and call it an
{\em augmented fan} where $\alpha$ is an imaginary cone.
We define $A(\alpha) := A$.
The category $\GEM(\hat\Delta )$ is defined similarly by supposing
$\sigma\prec\alpha$ for all $\sigma\in\hat\Delta$.
An object $L$ of $\GEM(\Delta)$ is also regarded as that of
$\GEM(\hat\Delta )$ by setting $L(\alpha) := \{0\}$.

For each $\rho\in\hat\Delta$, an additive covariant functor
\begin{equation}
\i_\rho^*\cl\GEM(\hat\Delta)\lra\GM(A(\rho))
\end{equation}
is defined by
\begin{equation}
\i_\rho^*(L) :=\bigoplus_{\sigma\in F(\rho)}%
L(\sigma)_{A(\rho)}\;,
\end{equation}
where $F(\rho)$ is the set of faces of $\rho$ and
we suppose $F(\rho) =\hat\Delta$ if $\rho =\alpha$.
Recall that $L(\sigma)_{A(\rho)}$ is the graded
$A(\rho)$-module $L(\sigma)\otimes_{A(\sigma)}A(\rho)$
for each $\sigma$.
We usually denote by $\Gamma$ the functor $\i_\alpha^*$

For a homomorphism $f\cl L\ra K$ in $\GEM(\hat\Delta)$, the
$(\sigma,\tau)$-component of the homomorphism
\begin{equation}
\begin{array}{ccccc}
\i_\rho^*(f) & \cl & \i_\rho^*(L) & \lra & \i_\rho^*(K) \\
            &     &     \|      &      &     \|      \\

 & & \displaystyle\bigoplus_{\sigma\in F(\rho)}L(\sigma)_{A(\rho)}
 & & \displaystyle\bigoplus_{\tau\in F(\rho)}K(\tau)_{A(\rho)}
\end{array}
\end{equation}
is defined to be $f(\sigma/\tau)_{A(\rho)}$ if
$\sigma\prec\tau$ and zero otherwise.

For each $\rho\in\Delta$, the additive covariant functor
\begin{equation}
\i_\rho^!\cl\GEM(\hat\Delta)\lra\GM(A(\rho))
\end{equation}
is defined by $\i_\rho^!(L) := L(\rho)$.
For $f\cl L\ra K$, the homomorphism $\i_\rho^!(f)$ is defined
to be $f(\rho/\rho)$.

For cones $\sigma,\tau$ with $\sigma\prec\tau$ and
$r_\tau = r_\sigma + 1$, we define the {\em incidence isomorphism}
$q'_{\sigma/\tau}\cl\det(\sigma)\ra\det(\tau)$ of free
$\Z$-modules of rank one as follows.

By the condition, $N(\tau)/N(\sigma)$ is a free $\Z$-module
of rank one.
We take $a\in N(\tau)\cap\tau$ such that the class of
$a$ in $N(\tau)/N(\sigma)$ is a generator.
Then we define $q'_{\sigma/\tau}(w) := a\wedge w$ for
$w\in\det(\sigma)$.

For cones $\sigma,\rho$ with $\sigma\prec\rho$ and
$r_\rho = r_\sigma + 2$, there exists exactly two cones
$\tau$ with $\sigma\prec\tau\prec\rho$ and
$r_\tau = r_\sigma + 1$.
Let these cones be $\tau_1,\tau_2$.
Then the equality
\begin{equation}
\label{codim 2 eq}
q'_{\sigma/\tau_1}\cdot q'_{\tau_1/\rho} +
q'_{\sigma/\tau_2}\cdot q'_{\tau_2/\rho} = 0
\end{equation}
holds (cf. \cite[Lem.1.4]{Ishida1}).

For a subset $\Phi\subset\Delta$ and an integer $i$, we set
\begin{equation}
\Phi(i) :=\{\sigma\in\Phi\scl r_\sigma = i\}\;.
\end{equation}

A subset $\Phi$ of $\Delta$ is said to be {\em locally star closed}
if $\sigma,\rho\in\Phi$, $\tau\in\Delta$ and $\sigma\prec\tau\prec\rho$
imply $\tau\in\Phi$.

For a locally star closed subset $\Phi$ of $\Delta$, we define
a complex $E(\Phi,\Z)^\bullet$ of free $\Z$-modules as follows.

For each integer $i$, we set
\begin{equation}
E(\Phi,\Z)^i :=\bigoplus_{\sigma\in\Phi(i)}\det(\sigma)\;.
\end{equation}
For $\sigma\in\Phi(i)$ and $\tau\in\Phi(i+1)$, the
$(\sigma,\tau)$-component of the coboundary map
\begin{equation}
d^i\cl E(\Phi,\Z)^i\ra E(\Phi,\Z)^{i+1}
\end{equation}
is defined to be $q'_{\sigma/\tau}$.
The equality $d^{i+1}\cdot d^i = 0$ follows from
(\ref{codim 2 eq}) for every $i$.

We say that a locally star closed subset $\Phi\subset\Delta$
is {\em $1$-complete} if, for each $\sigma\in\Phi(r-1)$, there exist
exactly two $\tau$'s in $\Phi(r)$ with $\sigma\prec\tau$.

If the finite fan $\Delta$ is complete \cite[Thm.1.11]{Oda1}, then
$\Delta(\sigma{\prec}) :=\{\rho\in\Delta\scl\sigma\prec\rho\}$
is $1$-complete for every $\sigma\in\Delta$.

If $\Phi$ is $1$-complete, then we can define an augmented complex
$E(\hat\Phi,\Z)^\bullet$ for $\hat\Phi :=\Phi\cup\{\alpha\}$ by
defining $r_\alpha := r+1$, $\det(\alpha) :=\bigwedge^r N$ and
$q'_{\tau/\alpha} := {\rm id}$ for every $\tau\in\Phi(r)$
with respect to the identification $\det(\tau) =\det(\alpha)$.
In particular, $E(\hat\Phi,\Z)^i = E(\Phi,\Z)^i$ for $i\not= r+1$
and $E(\hat\Phi,\Z)^{r+1} =\det(\alpha)$.

When $\Delta$ is complete, $E(\hat\Delta(\sigma{\prec}),\Z)^\bullet$
is acyclic for every $\sigma\in\Delta$.
Actually, $\H^i(E(\Delta(\sigma{\prec}),\Z)^\bullet)$ is equal
to the $(i{-}r_\sigma{-}1)$-th reduced cohomology group
of an $(r{-}r_\sigma{-}1)$-dimensional sphere, and hence it
vanishes if $i\not= r$.
The $r$-th cohomology is killed by
$E(\hat\Phi,\Z)^{r+1} =\det(\alpha)$.

We denote by $\CGEM(\Delta)$ and $\CGEM(\hat\Delta)$ the category
of finite $d$-complexes in $\GEM(\Delta)$ and $\GEM(\hat\Delta)$,
respectively.

Let $(L^\bullet, d_L)$ be an object of $\CGEM(\Delta)$.
Then, for each $\rho\in\Delta$, we get an object
$(L(\rho)^\bullet, d_L(\rho/\rho))$ of $\CGM(A(\rho))$
which we denote simply $L(\rho)^\bullet$.

For $\rho,\mu\in\Delta$ with $\rho\prec\mu$, we set
$F[\rho,\mu] :=\{\sigma\in F(\mu)\scl\rho\prec\sigma\}$.
Then the equality $d_L\cdot d_L = 0$ implies that
\begin{equation}
\label{CGEM1}
\sum_{\sigma\in F[\rho,\mu]}%
d_L^{i+1}(\sigma/\mu)\cdot d_L^i(\rho/\sigma) = 0
\end{equation}
for each integer $i$.
In particular, if $r_\mu =r_\rho + 1$, then the collection
$(d_L^i(\rho/\mu)\scl i\in\Z)$ defines a homomorphism of
complexes $d_L(\rho/\mu)\cl L(\rho)^\bullet\ra L(\mu)[1]^\bullet$
since then $F[\rho,\mu] =\{\rho,\mu\}$ and the equality (\ref{CGEM1})
imply the commutativity of the diagram
\begin{equation}
\begin{array}{ccc}
   \makebox[40pt]{}L(\rho)^i
 & \mathop{\lra}\limits^{\textstyle d_L^i(\rho/\rho)}
 & L^{i+1}(\rho)\makebox[70pt]{} \\
    d_L^i(\rho/\mu)\downarrow  &
 & \downarrow d_L^{i+1}(\rho/\mu)\\
   \makebox[40pt]{}L(\mu)[1]^i
 &\mathop{\lra}\limits^{\textstyle d_{L[1]}^i(\mu/\mu)}
 & L(\mu)[1]^{i+1}\makebox[70pt]{}
\end{array}
\;,
\end{equation}
where $d_{L[1]}^i(\mu/\mu) = - d_{L}^{i+1}(\mu/\mu)$.

Conversely, assume that complexes
$L(\rho)^\bullet\in\CGM(A(\rho))$ for $\rho\in\Delta$
and homomorphisms
\begin{equation}
d_L^i(\sigma/\tau)\cl L(\sigma)^i\lra L(\tau)^{i+1}
\end{equation}
for $\sigma,\tau\in\Delta$ with $\sigma\prec\tau$ and
$i\in\Z$ are given.
If they satisfy (\ref{CGEM1}) for all $(\rho,\mu)$ and
$i\in\Z$, then we get a complex $(L^\bullet, d_L)$ in
$\CGEM(\Delta)$.

An object $L^\bullet$ in $\CGEM(\Delta)$ is said to be
{\em shallow} if $d_L(\sigma/\rho) = 0$ for any $\sigma,\rho$
with $r_\rho - r_\sigma\geq 2$.

In order to define a shallow object $L^\bullet$ of
$\CGEM(\Delta)$, it is sufficient to give the following data
(1), (2) and check the condition (3).

(1) A complex $L(\sigma)^\bullet\in\CGM(A(\sigma))$ for
each $\sigma\in\Delta$.

(2) A homomorphism
$d(\sigma/\tau)\cl L(\sigma)^\bullet\ra L(\tau)[1]^\bullet$
for each pair $(\sigma,\tau)$ of cones in $\Delta$ with
$\sigma\prec\tau$ and $r_\tau = r_\sigma + 1$.

(3) The equality
\begin{equation}
d(\tau_1/\rho)^{i+1}\cdot d(\sigma/\tau_1)^i +
d(\tau_2/\rho)^{i+1}\cdot d(\sigma/\tau_2)^i = 0
\end{equation}
holds for all $i\in\Z$ and all pairs $(\sigma,\rho)$ with
$\sigma\prec\rho$ and $r_\rho = r_\sigma + 2$, where
$\tau_1,\tau_2$ are the dual cones with
$\sigma\prec\tau_i\prec\rho$ and $r_{\tau_i} = r_\sigma + 1$.

For an object $L^\bullet$ in $\CGEM(\Delta)$, we define
a shallow object $\tilde L^\bullet$ as follows.

For each $\rho\in\Delta$ and $i\in\Z$, we set
\begin{equation}
\label{def of tilde L}
\tilde L(\rho)^i :=
\bigoplus_{\sigma\in F(\rho)}\;\bigoplus_{\eta\in F(\sigma)}
\det(\rho)\otimes\det(\sigma)^*\otimes
L(\eta)_{A(\rho)}^{-r_\rho+r_\sigma+i}\;,
\end{equation}
where $\det(\sigma)^* :=\Hom_\Z(\det(\sigma),\Z)$.
For
\begin{equation}
\label{def of tilde L +1}
\tilde L(\rho)^{i+1} :=
\bigoplus_{\tau\in F(\rho)}\;\bigoplus_{\zeta\in F(\tau)}
\det(\rho)\otimes\det(\tau)^*\otimes
L(\zeta)_{A(\rho)}^{-r_\rho+r_\tau+i+1}\;,
\end{equation}
the $((\sigma,\eta),(\tau,\zeta))$-component of the
coboundary map
$d(\rho/\rho)^i\cl \tilde L(\rho)^i\ra \tilde L(\rho)^{i+1}$
is defined to be zero map except for the cases (a) $\sigma =\tau$ and
$\eta\prec\zeta$, or (b) $\tau\prec\sigma$, $r_\sigma = r_\tau + 1$
and $\eta =\zeta$. In case (a), the component is defined to be
$(-1)^{r_\rho - r_\sigma}{\rm id}
\otimes d(\eta/\zeta)^{-r_\rho+r_\sigma+i}$, and
in case (b), it is defined to be $(-1)^{r_\rho - r_\sigma-1}
1_{\det(\rho)}\otimes(q'_{\tau/\sigma})^*\otimes{\rm id}$, where
${\rm id}$'s are identity maps of the corresponding parts.

For $\rho\prec\mu$ with $r_\mu = r_\rho + 1$, the homomorphism
$d(\rho/\mu)\cl\tilde L(\rho)^\bullet\ra\tilde L(\mu)[1]^\bullet$
is defined to be the tensor product of $q'_{\rho/\mu}$ and the
natural inclusion map.

It is easy to check that $\tilde L(\rho)^\bullet$ is actually
a complex and $d(\rho/\mu)$ is a homomorphism of complexes.
The condition (3) in the definition of shallow complexes is
satisfied by the equality (\ref{codim 2 eq}).

We define a homomorphism $f_L\cl L^\bullet\ra\tilde L^\bullet$
as follows.

For $\tau,\rho\in\Delta$ with $\tau\prec\rho$, the
$A(\tau)$-homomorphism
$f_L(\tau/\rho)^i\cl L(\tau)^i\ra\tilde L(\rho)^i$ is the
inclusion map to the component
$\det(\rho)\otimes\det(\sigma)^*\otimes
L(\eta)_{A(\rho)}^{-r_\rho+r_\sigma+i}$ for $\sigma =\rho$
and $\eta =\tau$ in the description (\ref{def of tilde L}).
The compatibility with the coboundary maps is checked easily.

A homomorphism $f\cl L^\bullet\ra K^\bullet$ of $d$-complexes of
graded exterior modules on $\Delta$ is said to be a
{\em quasi-isomorphism} if
$f(\sigma/\sigma)\cl L(\sigma)^\bullet\ra K(\sigma)^\bullet$ is a
quasi-isomorphism of $d$-complexes of $A(\sigma)$-modules
for every $\sigma\in\Delta$.

The following proposition shows that any object in $\CGEM(\Delta)$
is quasi-isomorphic to a shallow one.

%
%

\begin{Prop}
\label{Prop 2.1}
For any $L^\bullet$ in $\CGEM(\Delta)$, the homomorphism
$f_L\cl L^\bullet\ra\tilde L^\bullet$ is quasi-isomorphic.
\end{Prop}

\Proof
We have to check that
$f_L(\rho/\rho)\cl L(\rho)^\bullet\ra\tilde L(\rho)^\bullet$
is a quasi-isomorphism for each $\rho\in\Delta$.
We define a decreasing filtration $\{F^k\}$ on the complex
$V^\bullet :=\tilde L(\rho)^\bullet$ as follows.
For each integer $k$, $F^k(V)^i$ is defined to be the direct
sum of the components
\begin{equation}
\det(\rho)\otimes\det(\sigma)^*\otimes
L(\eta)_{A(\rho)}^{-r_\rho+r_\sigma+i}
\end{equation}
in the description (\ref{def of tilde L}) for $\sigma,\eta$
with $r_\eta\geq k$.
Then $F^k(V)^0 = V^\bullet$ and $F^k(V)^\bullet$ is a
subcomplex of $V^\bullet$ for each $k\geq 0$.
We can decompose the quotient complex
$F^k(V)^\bullet/F^{k+1}(V)^\bullet$ to a direct sum
\begin{equation}
\bigoplus_{\eta\in F(\rho)(k)}V_{k,\eta}^\bullet\;,
\end{equation}
where $V_{k,\eta}^\bullet$ is the part consisting of the
components related to $\eta$.
Then we see that $V_{k,\eta}^\bullet$ is isomorphic to the
associated single complex of the bicomplex
\begin{equation}
\det(\rho)\otimes\Hom(E(F[\eta,\rho],\Z),\Z)[-r_\rho]^\bullet
\otimes L(\eta)_{A(\rho)}^\bullet\;,
\end{equation}
where we denote by $\Hom(E,\Z)^\bullet$ the complex defined by
\begin{equation}
\Hom(E,\Z)^i :=\Hom(E^{-i},\Z)
\end{equation}
and $d^i := (-1)^{i+1}(d_E^{-i-1})^*$ for $i\in\Z$.
If $k < r_\rho$, then $\eta\not=\rho$ and $E(F[\eta,\rho],\Z)^\bullet$
is acyclic, and hence so is $V_{k,\eta}^\bullet$.
Hence $V^\bullet$ is quasi-isomorphic to the subcomplex
$F^{r_\rho}(V)^\bullet$.
Since $f_L(\rho/\rho)$ is an isomorphism from $L(\rho)^\bullet$
onto $F^{r_\rho}(V)^\bullet$, it is a quasi-isomorphism to
$V^\bullet =\tilde L(\rho)^\bullet$.
\QED

For a complex $L^\bullet$ in $\CGEM(\Delta)$, we define a
shallow complex $\D(L)^\bullet$ in $\CGEM(\Delta)$ as follows.

For each $\rho\in\Delta$, we set
\begin{equation}
\D(L)(\rho)^\bullet :=
\det(\rho)\otimes\d_\rho(\i_\rho^*L)[-r_\rho]^\bullet\;.
\end{equation}

Let $\rho,\mu\in\Delta$ satisfy $\rho\prec\mu$ and
$r_\mu = r_\rho + 1$.
By the definition of $\i_\mu^*$ and $\i_\rho^*$, we see that
$(\i_\rho^*L^i)_{A(\mu)}$ is a direct summand of
$\i_\mu^*L^i$.
Let $p(\mu/\rho)^i\cl\i_\mu^*L^i\ra(\i_\rho^*L^i)_{A(\mu)}$ be
the natural projection map.
We see that $\{p(\mu/\rho)^i\scl i\in\Z\}$ defines a homomorphism
$p(\mu/\rho)\cl\i_\mu^*L^\bullet\ra(\i_\rho^*L)_{A(\mu)}^\bullet$.
Hence we get a homomorphism
\begin{equation}
\label{dual rho mu}
\d_\mu(p(\mu/\rho))\cl\d_\mu((\i_\rho^*L)_{A(\mu)})^\bullet\lra
\d_\mu(\i_\mu^*L)^\bullet\;.
\end{equation}
Since
\begin{equation}
\d_\mu((\i_\rho^*L)_{A(\mu)})^\bullet =
\d_\rho(\i_\rho^*L)_{A(\mu)}^\bullet
\end{equation}
by Lemma~\ref{Lem 1.4}, we get a homomorphism
\begin{equation}
i(\rho/\mu)\cl\d_\rho(\i_\rho^*L)^\bullet\lra
\d_\mu(\i_\mu^*L)^\bullet
\end{equation}
as the composite of the inclusion
$\d_\rho(\i_\rho^*L)^\bullet\ra\d_\rho(\i_\rho^*L)_{A(\mu)}^\bullet$
and (\ref{dual rho mu}).
We define
\begin{equation}
d_{\D(L)}(\rho/\mu)\cl\D(L)(\rho)^\bullet\lra
\D(L)(\mu)[1]^\bullet
\end{equation}
to be $q_{\rho/\mu}\otimes i(\rho/\mu)[-r_\rho]$.
By the equality (\ref{codim 2 eq}), the condition (3)
of the construction of shallow complexes
is satisfied and we get a complex $\D(L)^\bullet$.

%
%

\begin{Lem}
\label{Lem 2.2}
Let $L^\bullet$ be an object of $\CGEM(\Delta)$.
Then, there exists a quasi-isomorphism
$\varphi\cl L^\bullet\ra \D(\D(L))^\bullet$.
\end{Lem}

\Proof
We prove that $\D(\D(L))^\bullet$ is isomorphic to
$\tilde L^\bullet$.
Then the lemma follows from Proposition~\ref{Prop 2.1}.

Let $\rho$ be an element of $\Delta$.
For each integer $i$,
\begin{eqnarray*}
 &   & \D(\D(L))(\rho)^i \\
 & = & \bigoplus_{\sigma\in F(\rho)}%
       \det(\rho)\otimes\d_\rho(\D(L)(\sigma)_{A(\rho)}^{r_\rho-i}) \\
 & = & \bigoplus_{\sigma\in F(\rho)}\bigoplus_{\eta\in F(\sigma)}%
       \det(\rho)\otimes\d_\rho(\det(\sigma)\otimes
       \d_\sigma(L(\eta)_{A(\sigma)}^{-r_\rho+r_\sigma+i})_{A(\rho)}) \\
 & = & \bigoplus_{\sigma\in F(\rho)}\bigoplus_{\eta\in F(\sigma)}%
       \det(\rho)\otimes\det(\sigma)^*\otimes\d_\rho(
       \d_\sigma(L(\eta)_{A(\sigma)}^{-r_\rho+r_\sigma+i})_{A(\rho)})\;.
\end{eqnarray*}
By Lemma~\ref{Lem 1.4}, we have a natural isomorphism
\begin{equation}
\d_\rho(\d_\sigma(L(\eta)_{A(\sigma)}^{-r_\rho+r_\sigma+i})_{A(\rho)})
\simeq\d_\rho(\d_\rho(L(\eta)_{A(\rho)}^{-r_\rho+r_\sigma+i}))\;.
\end{equation}
We identify the last $A(\rho)$-module with
$L(\eta)_{A(\rho)}^{-r_\rho+r_\sigma+i}$ by the canonical
isomorphism for all $(\sigma,\eta)$.
we know $\D(\D(L))(\rho)^i =\tilde L(\rho)^i$ for every
$\rho\in\Delta$ and $i\in\Z$, however the coboundary map is not
equal to that of $\tilde L^\bullet$.

We consider the descriptions (\ref{def of tilde L}) and
(\ref{def of tilde L +1}) with replacing $\tilde L(\rho)^i$ and
$\tilde L(\rho)^{i+1}$ by $\D(\D(L))(\rho)^i$ and
$\D(\D(L))(\rho)^{i+1}$, respectively.
We can check that the $((\sigma,\eta),(\tau,\zeta))$-component of the
coboundary map
$d(\rho/\rho)^i\cl \tilde \D(\D(L))(\rho)^i\ra\D(\D(L))(\rho)^{i+1}$
is the zero map except for the cases (a) $\sigma =\tau$ and
$\eta\prec\zeta$, or (b) $\tau\prec\sigma$, $r_\sigma = r_\tau + 1$
and $\eta =\zeta$. In case (a), the component is calculated to be
$(-1)^{r_\rho + 1}{\rm id}\otimes d(\eta/\zeta)^{-r_\sigma+r_\sigma+i}$,
and in case (b), it is
$(-1)^{i+1}1_{\det(\rho)}\otimes(q'_{\tau/\sigma})^*\otimes{\rm id}$.

For each $\rho$ and $i$, we define an isomorphism
\begin{equation}
\varphi(\rho)^i\cl\tilde L(\rho)^i\lra\D(\D(L))(\rho)^i
\end{equation}
by defining its restriction to the component
\begin{equation}
\det(\rho)\otimes\det(\sigma)^*\otimes
L(\eta)_{A(\rho)}^{-r_\rho+r_\sigma+i}
\end{equation}
to be $(-1)^{(r_\rho+i)(r_\sigma+1)}$ times the
identity map to the same component of $\D(\D(L))(\rho)^i$.
Then we see that the collection $\{\varphi(\rho)^i\}$ defines
an unmixed isomorphism $\tilde L^\bullet\ra\D(\D(L))^\bullet$.
\QED

We can define a similar functor
\begin{equation}
\hat\D\cl\CGEM(\hat\Delta)\lra\CGEM(\hat\Delta)
\end{equation}
for an augmented 1-complete fan $\hat\Delta$ by the convention
that $r_\alpha = r+1$, $F(\alpha) =\hat\Delta$ and
$q_{\tau/\alpha} := 1_{\det N}$ for $\tau\in\Delta(r)$.
When $L^\bullet$ is an object of $\CGEM(\Delta)$,
$\hat\D(L)^\bullet$ is in $\CGEM(\hat\Delta)$.
Since $\det(\alpha) =\det N$, we get an exact sequence
\begin{equation}
0\ra\hat\D(L)(\alpha)^\bullet\lra
\hat\D(L)^\bullet\lra\D(L)^\bullet\ra 0
\end{equation}
in $\CGEM(\hat\Delta)$ as well as an exact sequence
\begin{equation}
\label{augmented exact}
0\ra\det N\otimes\d_N(\Gamma(L))[-r-1]^\bullet\lra
\Gamma(\hat\D(L))^\bullet\lra\Gamma(\D(L))^\bullet\ra 0
\end{equation}
in $\CGM(A)$ by applying $\Gamma$.

%

\begin{Lem}
\label{Lem 2.3}
Assume that $\Delta$ is a complete fan.
Then, for any $L^\bullet$ in $\CGEM(\Delta)$, the complex
of $A$-modules $\Gamma((\hat\D(L)))^\bullet$ is acyclic.
\end{Lem}

\Proof
For each $i\in\Z$, we have
\begin{equation}
\label{hat dual}
\Gamma(\hat\D(L))^i =
\bigoplus_{\rho\in\hat\Delta}\;\bigoplus_{\sigma\in F(\rho)}
\det(\rho)\otimes\d_\rho(L(\sigma)_{A(\rho)}^{-r_\rho-i})_A\;.
\end{equation}
For
\begin{equation}
\Gamma(\hat\D(L))^{i+1} =
\bigoplus_{\mu\in\hat\Delta}\;\bigoplus_{\tau\in F(\mu)}
\det(\mu)\otimes\d_\mu(L(\tau)_{A(\mu)}^{-r_\mu-i-1})_A\;,
\end{equation}
the $((\rho,\sigma),(\mu,\tau))$-component of the coboundary
map is nonzero only for (a) $\rho =\mu$ and $\tau\prec\sigma$,
or (b) $\rho\prec\mu$, $r_\mu = r_\rho + 1$ and $\sigma =\tau$.
In case (a), the component is
$(-1)^{i+1}{\rm id}\otimes\d_\rho(d(\tau/\sigma)_{A(\rho)}^{-r_\rho-i-1})_A$,
and in case (b), it is $q'_{\rho/\mu}\otimes{\rm id}$, where
${\rm id}$'s are the identities of the corresponding parts,
respectively.

For $V^\bullet :=\Gamma(\hat\D(L))^\bullet$, we introduce a
deceasing filtration $\{F^k\}$ as follow.
For each integer $i$, we define $F^k(V)^i$ to be the direct
sum of the components
$\det(\rho)\otimes\d_\rho(L(\sigma)_{A(\rho)}^{-r_\rho-i})_A$ of
$\Gamma(\hat\D(L))^i$ in the description (\ref{hat dual})
for all pairs $(\rho,\sigma)$ with $r_\sigma\geq k$.
Then $F^k(V)^\bullet$ is a subcomplex of $V^\bullet$ for
every $k\in\Z$, and $F^k(V)^\bullet/F^{k+1}(V)^\bullet$ is
a direct sum of complexes
$\bigoplus_{\sigma\in\Delta(k)}V_\sigma^\bullet$, where
$V_\sigma^\bullet$ is the part related to each $\sigma$.
We see that $V_\sigma^\bullet$ is isomorphic to the
associated single complex of the bicomplex
\begin{equation}
E(F[\sigma,\alpha],\Z)^\bullet\otimes
\d_\sigma(L(\sigma))_A^\bullet\;.
\end{equation}
Since $\Delta$ is complete, $E(F[\sigma,\alpha],\Z)^\bullet$
is acyclic for every $\sigma\in\Delta$.
Hence $F^k(V)^\bullet/F^{k+1}(V)^\bullet$ is also acyclic
for every $k$.
Since $F^0(V)^\bullet = V^\bullet$ and $F^{r+1}(V)^\bullet =\{0\}$,
$V^\bullet$ is acyclic.
\QED

By the long exact sequence obtained from the exact sequence
(\ref{augmented exact}), we get the following corollary.

%
%

\begin{Cor}
\label{Cor 2.4}
If $\Delta$ is complete, then there exists an isomorphism
\begin{equation}
\H^{-p}(\d_N(\Gamma(L))^\bullet)\simeq
\H^{r-p}(\Gamma(\D(L))^\bullet)
\end{equation}
in $\GM(A)$ for each integer $p$.
\end{Cor}

The following proposition is a consequence of Lemma~\ref{Lem 1.5}
and Corollary~\ref{Cor 2.4}.

%
%

\begin{Prop}
\label{Prop 2.5}
Assume that $\Delta$ is complete.
For any integer $p, q$, the equality
\begin{equation}
\dim_\Q\H^p(\Gamma(L)^\bullet)_q =
\dim_\Q\H^{r-p}(\Gamma(\D(L))^\bullet)_{-r-q}
\end{equation}
holds.
\end{Prop}

A finite fan $\Delta$ of $N_\R$ is said to be {\em lifted complete}
if there exists a rational line $\ell$ of $N_\R$ going through
the origin with the following property.

Let $\bar\sigma$ be the image of $\sigma$ in the quotient $N_\R/\ell$
for each $\sigma\in\Delta$.
Then (1) $\dim\bar\sigma =r_\sigma$ for every $\sigma\in\Delta$,
(2) $\bar\sigma\not=\bar\tau$ for any distinct $\sigma,\tau\in\Delta$
and (3) $\Delta_\ell :=\{\bar\sigma\scl\sigma\in\Delta\}$ is a
complete fan of $N_\R/\ell$.

For a lifted complete fan, the associated toric variety has
an action of the multiplicative algebraic group $\Gm$, and it
has a complete toric variety of dimension $r-1$ as the
geometric quotient in the sense of Mumford's geometric invariant
theory.

Let $\Delta$ be a lifted complete fan with respect to $\ell$ and
let $\ell^+$ be one of the one-dimensional cones contained
in $\ell$.
By the property (1), $\dim(\tau +\ell^+) =r_\tau + 1$ for every
$\tau\in\Delta$.
The {\em oriented} lifted complete fan $\tilde\Delta$ is defined
to be $\Delta\cup\{\beta\}$ where $\beta$ is an imaginary
cone of dimension $r$.
We suppose $\sigma\prec\beta$ for every $\sigma\in\Delta$.
We define $\det(\beta) :=\det N$ and
$q'_{\tau/\beta} := q'_{\tau/\tau'}$ for
$\tau\in\Delta(r{-}1)$, where $\tau' :=\tau +\ell^+$

For any star closed subset $\Phi\subset\tilde\Delta$, the
complex $E(\Phi,\Z)^\bullet$ is defined similarly as in the
previous case.
The complex $E(\tilde\Delta(\sigma{\prec}),\Z)^\bullet$ is
acyclic for every $\sigma\in\Delta$ since it is isomorphic to
the augmented complex
$E(\hat\Delta_\ell,\Z)^\bullet$.

Each lifted complete fan has two orientations according to the
choice of $\ell^+$, but it does not depend on the choice
of the line $\ell$.

Here we give three typical examples of oriented lifted complete fans.

(1)
Let $C\subset N_\R$ be a closed convex cone of dimension $r$
which may not be strongly convex and which is not equal to $N_\R$.
Let $\partial C$ be the boundary set of $C$.
Then a finite fan $\Delta$ with the support $\partial C$ is
lifted complete.
Any $\ell$ which intersects the interior of $C$ satisfies the
condition.
As the natural orientation, we take $\ell^+ :=\ell\cap C$.

(2)
Let $\Phi$ be a simplicial complete fan of $N_\R$ and let $\gamma$
be a one-dimensional cone in $\Phi$.
Set
\begin{equation}
\Delta :=
\{\sigma\in\Phi\scl\gamma\not\prec\sigma,\sigma+\gamma\in\Phi\}\;.
\end{equation}
Then $\Delta$ is a lifted complete fan and $\ell^+ :=\gamma$
defines an orientation (cf.\cite{Oda2}).

(3)
Let $N'$ be a free $\Z$-module of rank $r - 1$, $\Phi$ a
finite complete fan of $N'_\R$ and $h$ a real-valued
continuous function on $N'_\R$ which is linear on each cone
$\sigma\in\Phi$ and has rational values on $N'_\Q$.
Then the fan $\Delta =\{\sigma'\scl\sigma\in\Phi\}$ of
$N_\R := N'_\R\oplus \R$ is lifted complete, where
$\sigma' :=\{(x,h(x))\scl x\in\sigma\}$ for each $\sigma\in\Phi$.
We take $\ell^+ :=\{0\}\times\R_0$ as the orientation.
This type of fan is treated in \cite{Oda3}.

Let $\tilde\Delta =\Delta\cup\{\beta\}$ be an oriented
lifted complete fan.
The category $\GEM(\tilde\Delta)$ is defined by setting
$A(\beta) := A$.
Then the functors
\begin{equation}
\i_\beta^*\cl\GEM(\tilde\Delta)\lra\GM(A)
\end{equation}
and
\begin{equation}
\tilde\D\cl\GEM(\tilde\Delta)\lra\GEM(\tilde\Delta)
\end{equation}
are defined similarly as $\i_\alpha^*$ and $\hat\D$ in
the case of augmented 1-complete fans, respectively.
We denote also by $\Gamma$ the functor $\i_\beta^*$.

We omit the proofs of the following results, since they are
similar to those of the corresponding results for an
augmented complete fan.

%
%

\begin{Lem}
\label{Lem 2.6}
Let $\tilde\Delta =\Delta\cup\{\beta\}$ be an oriented
lifted complete fan.
Then, for any $L^\bullet$ in $\CGEM(\Delta)$, the complex
of $A$-modules $\Gamma(\tilde\D(L))^\bullet$ is acyclic, and
there exists an exact sequence
\begin{equation}
0\lra\det N\otimes\d(\Gamma(L))^\bullet[-r]\lra
\Gamma(\tilde\D(L))^\bullet\lra
\Gamma(\D(L))^\bullet\lra 0\;.
\end{equation}
\end{Lem}

%
%

\begin{Cor}
\label{Cor 2.7}
Let $\Delta\cup\{\beta\}$ be an oriented lifted complete fan.
Then for any $L^\bullet\in\CGEM(\Delta)$, there exists an isomorphism
\begin{equation}
\H^{-p}(\d_N(\Gamma(L))^\bullet)\simeq
\H^{r-p-1}(\Gamma(\D(L))^\bullet)
\end{equation}
in $\GM(A)$ for each integer $p$.
\end{Cor}

%
%

\begin{Prop}
\label{Prop 2.8}
Let $\Delta\cup\{\beta\}$ be an oriented lifted complete fan.
For any integer $p, q$, the equality
\begin{equation}
\dim\H^p(\Gamma(L)^\bullet)_q =
\dim\H^{r-p-1}(\Gamma(\D(L))^\bullet)_{-r-q}
\end{equation}
holds.
\end{Prop}

Let $\Delta$ be a finite fan of $N_\R$ and let $\Phi$ be
a subfan of $\Delta$.
For $K^\bullet$ in $\CGEM(\Phi)$, we denote by the same symbol
$K^\bullet$ the trivial extension to $\Delta$, i.e., we define
$K(\sigma)^\bullet :=\{0\}$ for $\sigma\in\Delta\setminus\Phi$.

For $L^\bullet\in\CGEM(\Delta)$ and $\sigma\in\Delta$, we
set
\begin{equation}
\i_\sigma^\circ(L)^\bullet :=
\i_\sigma^*(L|F(\sigma)\setmin\{\sigma\})^\bullet\;,
\end{equation}
where $(L|F(\sigma)\setmin\{\sigma\})^\bullet$ is the
restriction of $L^\bullet$ to $F(\sigma)\setminus\{\sigma\}$.

We see the case $\Phi =\Delta\setminus\{\pi\}$ for a maximal element
$\pi\in\Delta$.
Let $L^\bullet\in\CGEM(\Delta)$.
Since
\begin{equation}
\i_\pi^\circ(L)^i =
\bigoplus_{\sigma\in F(\pi)\setminus\{\pi\}}L(\sigma)_{A(\pi)}^i
\end{equation}
and $\i_\pi^*(L)^i = L(\pi)^i\oplus\i_\pi^\circ(L)^i$ for
each $i\in\Z$, there exists an exact sequence
\begin{equation}
0\ra L(\pi)^\bullet\lra\i_\pi^*(L)^\bullet\lra
\i_\pi^\circ(L)^\bullet\ra 0
\end{equation}
in $\CGM(A(\pi))$.
In other words, $\i_\pi^*(L)^\bullet$ is equal to the
mapping cone of the homomorphism
\begin{equation}
\label{hom of ext}
\phi\cl\i_\pi^\circ(L)[-1]^\bullet\lra L(\pi)^\bullet
\end{equation}
whose component for each $\sigma\in F(\pi)\setminus\{\pi\}$ is
$d(\sigma/\pi)_{A(\pi)}$.

The extension of $(L|\Phi)^\bullet$ to $L^\bullet$ is determined
by the above homomorphism $\phi$.
Actually, if $(L|\Phi)^\bullet\in\CGEM(\Phi)$,
$L(\pi)^\bullet$ and the homomorphism (\ref{hom of ext}) is
given, then we get the extension $L^\bullet$.

We define a functor $\j_!^\Phi\cl\CGEM(\Phi)\ra\CGEM(\Delta)$
as follows.

For $K^\bullet$ in $\CGEM(\Phi)$, we define $\j_!^\Phi(K)^\bullet$
by
\begin{equation}
\j_!^\Phi(K)(\sigma)^\bullet := K(\sigma)^\bullet
\end{equation}
for $\sigma\in\Phi$ and
\begin{equation}
\j_!^\Phi(K)(\pi)^\bullet := \i_\pi^\circ(K)[-1]^\bullet
= \i_\pi^*(K)[-1]^\bullet\;.
\end{equation}
In particular,
\begin{equation}
\j_!^\Phi(K)(\pi)^{i+1} =
\bigoplus_{\sigma\in F(\pi)\setminus\{\pi\}}K(\sigma)_{A(\pi)}^i\;.
\end{equation}
We define the extension $\j_!^\Phi(K)^\bullet\in\CGEM(\Delta)$
of $K^\bullet$ by the identity map
\begin{equation}
\i_\pi^\circ(K)[-1]^\bullet\lra\j_!^\Phi(K)(\pi)^\bullet\;.
\end{equation}

Let $(L^\bullet, d_L)$ and $(K^\bullet, d_K)$ be $d$-complexes
in $\CGEM(\Delta)$.
By an {\em unmixed} homomorphsm $f\cl L^\bullet\ra K^\bullet$, we means a
collection $\{f^i\scl i\in\Z\}$ of unmixed homomorphisms $f^i\cl L^i\ra K^i$
such that $d_K^i\cdot f^i = f^{i+1}\cdot d_L^i$ for every $i\in\Z$.
Note that $d_L^i$ and $d_K^i$ are not necessary unmixed.

Let $L^\bullet$ be an object of $\CGEM(\Delta)$.
Then there exists a unique unmixed homomorphism
$\j_!^\Phi(L|\Phi)^\bullet\ra L^\bullet$ which is the
identity map on $\Phi$.
The homomorphism
\begin{equation}
\j_!^\Phi(L|\Phi)(\pi)^\bullet =\i_\pi^\circ(L)[-1]^\bullet
\ra L(\pi)^\bullet
\end{equation}
is defined to be (\ref{hom of ext}).

We can define a functor $\j_!^\Phi\cl\CGEM(\Phi)\ra\CGEM(\Delta)$
for general subfan $\Phi\subset\Delta$ as the composite of
the above functors.
However, we will not use the general case.

Let $\pi$ be a maximal element of $\Delta$.
For $L^\bullet$ in $\CGEM(\Delta)$ and an integer $k$, we define
$\grt_\pi^{\geq k}L^\bullet$ by
\begin{equation}
(\grt_\pi^{\geq k}L)(\sigma)^\bullet =
\left\{
    \begin{array}{ll}
        L(\sigma)^\bullet & \hbox{if } \sigma\not=\pi \\
        \grt^{\geq k}L(\pi)^\bullet & \hbox{if } \sigma=\pi\;.
    \end{array}
\right.
\end{equation}
There exists a natural unmixed homomorphism
$L^\bullet\lra\grt_\pi^{\geq k}L^\bullet$
such that the components for $\sigma\not=\pi$ are the identity
maps and the component for $\pi$ is the natural surjection
$L^\bullet(\pi)\ra\grt^{\geq k}(L(\pi))^\bullet$.
This homomorphism is a quasi-isomorphism if and only if
$\grt_{\leq k-1}(L^\bullet(\pi))$ is acyclic.

Two objects $L^\bullet, K^\bullet$ in $\CGEM(\Delta)$ are
said to be {\em quasi-isomorphic} if there exist a finite
sequence $L_0^\bullet, L_1^\bullet,\cdots, L_{2k}^\bullet$
of objects in $\CGEM(\Delta)$ with $L^\bullet = L_0^\bullet$
and $L_{2k}^\bullet = K^\bullet$ and quasi-isomorphisms
$L_{2i-2}^\bullet\ra L_{2i-1}^\bullet$ and
$L_{2i}^\bullet\ra L_{2i-1}^\bullet$ for $i = 1,\cdots, k$.
Some lemmas on this definition are given at the end this section.

A map $\p\cl\Delta\setminus\{\Zero\}\ra\Z$ is called
a {\em perversity} on a finite fan $\Delta$.
We prove the following theorem.

%
%

\begin{Thm}
\label{Thm 2.9}
Let $\p$ be a perversity on $\Delta$.
Then there exists a finite $d$-complex $\ic_\p(\Delta)^\bullet$
of graded exterior modules on $\Delta$ satisfying
the following conditions.

{\rm (1)} $\H^0(\ic_\p(\Delta)(\Zero)^\bullet) =\Q$ and
$\H^i(\ic_\p(\Delta)(\Zero)^\bullet) =\{0\}$ for $i\not= 0$.

{\rm (2)} For $\sigma\in\Delta\setminus\{\Zero\}$ and
$i, j\in\Z$ with $i + j\leq\p(\sigma)$, we have
$\H^i(\ic_\p(\Delta)(\sigma)^\bullet)_j =\{0\}$.

{\rm (3)} For $\sigma\in\Delta\setminus\{\Zero\}$ and
$i, j\in\Z$ with $i + j\geq\p(\sigma)$, we have
$\H^i(\i_\sigma^*(\ic_\p(\Delta))^\bullet)_j =\{0\}$.

Futhermore, if $L^\bullet$ is another finite $d$-complex
satisfying the above conditions, then $L^\bullet$ is
quasi-isomorphic to $\ic_\p(\Delta)^\bullet$.
\end{Thm}

\Proof
We prove the theorem by induction on the number of
cones in $\Delta$.
If $\Delta =\{\Zero\}$, then we set
$\ic_\p(\Delta)(\Zero)^0 =\Q$ and
$\ic_\p(\Delta)(\Zero)^i =\{0\}$ for $i\not= 0$.
Then (1) and the last assertion are clearly satisfied.

Assume that $\Delta\not=\{\Zero\}$.
Let $\pi\in\Delta$ be a cone of maximal dimension.
We assume that the $d$-complex $\ic_\p(\Phi)^\bullet$ exists
for $\Phi =\Delta\setminus\{\pi\}$

Let $L^\bullet :=\j_!^\Phi\ic_\p(\Phi)^\bullet$.
We define
\begin{equation}
\ic_\p(\Delta)^\bullet :=\grt_\pi^{\geq\p(\pi)+1}L^\bullet\;.
\end{equation}

Since
\begin{equation}
\ic_\p(\Delta)(\pi)^\bullet =
\grt^{\geq\p(\pi)+1}L(\pi)^\bullet\;,
\end{equation}
the truncation $\grt_{\leq\p(\pi)}(\ic_\p(\Delta)(\pi))^\bullet$
is the zero complex.

On the other hand, $\i_\pi^*\ic_\p(\Delta)^\bullet$ is equal to
the mapping cone
\begin{equation}
(L(\pi)[1]\oplus\grt^{\geq\p(\pi)+1}L(\pi))^\bullet
\end{equation}
of the natural surjection
\begin{equation}
L(\pi)^\bullet\lra\grt^{\geq\p(\pi)+1}L(\pi)^\bullet\;.
\end{equation}
There exists an exact sequence
\begin{equation}
0\lra\widetilde{\grt}_{\leq\p(\pi)}L(\pi)^\bullet\lra L(\pi)^\bullet
\lra\grt^{\geq\p(\pi)+1}L(\pi)^\bullet\lra 0
\end{equation}
of $d$-complexes in $\GM(A(\pi))$.
Hence $\i_\pi^*\ic_\p(\Delta)^\bullet$ is quasi-isomorphic to
\begin{equation}
(\widetilde{\grt}_{\leq\p(\pi)}L(\pi))[1]^\bullet =
\widetilde{\grt}_{\leq\p(\pi)-1}(L(\pi)[1])^\bullet
\end{equation}
and to $\grt_{\leq\p(\pi)-1}(L(\pi)[1])^\bullet$.
Hence $\grt^{\geq\p(\pi)}\i_\pi^*\ic_\p(\Delta)^\bullet$ is acyclic.
This is equivalent to the condition (3).

For the last assertion, it is sufficient to prove the following
lemma.

%
%

\begin{Lem}
\label{Lem 2.10}
Let $\pi\not=\Zero$ be a maximal element of $\Delta$ and let
$\Phi :=\Delta\setminus\{\pi\}$.
Let $L^\bullet, K^\bullet$ be objects of $\CGEM(\Delta)$ which
satisfy the conditions of {\rm Theorem~\ref{Thm 2.9}}, and
assume that there exists a quasi-isomorphism
$(L|\Phi)^\bullet\ra(K|\Phi)^\bullet$.
Then, $L^\bullet$ and $K^\bullet$ is connected by a finite
sequence of unmixed quasi-isomorphisms.
\end{Lem}

\Proof
By the condition (2), the natural homomorphism
\begin{equation}
L^\bullet\ra\grt_\pi^{\geq \p(\sigma)+1}L^\bullet
\end{equation}
is a quasi-isomorphism.
Since $\i_\pi^*(L)^\bullet$ is the mapping cone of the
homomorphism
$\j_!^\Phi(L|\Phi)^\bullet\ra L(\pi)^\bullet$,
the condition (3) implies the homomorphism
\begin{equation}
\grt_\pi^{\geq \p(\sigma)+1}(\j_!^\Phi(L|\Phi))^\bullet\ra
\grt_\pi^{\geq \p(\sigma)+1}L^\bullet
\end{equation}
is also a quasi-isomorphism.

There are similar quasi-isomorphisms for $K^\bullet$.
We are done since the quasi-isomorphism
$(L|\Phi)^\bullet\ra(K|\Phi)^\bullet$ induces a quasi-isomorphism
\begin{equation}
\grt_\pi^{\geq \p(\sigma)+1}(\j_!^\Phi(L|\Phi))^\bullet\ra
\grt_\pi^{\geq \p(\sigma)+1}(\j_!^\Phi(K|\Phi))^\bullet\;.
\end{equation}
\QED

Thus we complete the proof of Theorem~\ref{Thm 2.9}

In the rest of this paper, we denote by $\ic_\p(\Delta)^\bullet$
the $d$-complex constructed in the proof of the above theorem,
and we call it the intersection complex of the fan $\Delta$
with the perversity $\p$.
It does not depend on the choice of the order of the induction,
and is uniquely determined by $\Delta$ and $\p$.
If $\Phi$ is a subfan of $\Delta$, then the construction
implies that $\ic_\p(\Phi)^\bullet$ is isomorphic to the
restriction of $\ic_\p(\Delta)^\bullet$ to $\Phi$.
In particular,
$\ic_\p(\Delta)(\sigma)^\bullet =\ic_\p(F(\sigma))(\sigma)^\bullet$.
Hence the complex $\ic_\p(\Delta)(\sigma)^\bullet$ of
$A(\sigma)$-modules depends only on $\p$ and $\sigma$.

%
%

\begin{Prop}
\label{Prop 2.11}
Let $\sigma$ be a nonzero cone of $\Delta$ and let $\p$
be a perversity on $\Delta$.
Then (1) $\ic_\p(\Delta)(\sigma)_j^i =\{0\}$ unless
$(i, j)\in [1, r_\sigma]\times[-r_\sigma, 0]$,
(2) $\i_\sigma^\circ(\ic_\p(\Delta))_j^i =\{0\}$ unless
$(i, j)\in [0, r_\sigma-1]\times[-r_\sigma, 0]$ and
(3) $\i_\sigma^*(\ic_\p(\Delta))_j^i =\{0\}$ unless
$(i, j)\in [0, r_\sigma]\times[-r_\sigma, 0]$.
Furthermore, $\Gamma(\ic_\p(\Delta))_j^i =\{0\}$ unless
$(i, j)\in [0, r]\times[-r, 0]$.
\end{Prop}

\Proof
We prove the proposition by induction on $r_\sigma$.
Note that $\ic_\p(\Delta)(\Zero)_j^i =\{0\}$ unless
$(i, j) = (0, 0)$ by the construction of $\ic_\p(\Delta)^\bullet$.

Let $\eta,\sigma$ be cones in $\Delta$ with $\eta\prec\sigma$.
Recall that, if we take an $(r_\sigma - r_\eta)$-dimensional
linear subspace $H$ of $N(\sigma)_\Q$ such that
$N(\sigma)_\Q = N(\eta)_\Q\oplus H$, then
$V_{A(\sigma)} = V\otimes_\Q A(H)$ for $V\in\GM(A(\eta))$.
Hence, if $V_j =\{0\}$ unless $j\in [a, b]$ for integers
$a, b$ with $a\leq b$, then $(V_{A(\sigma)})_j =\{0\}$
unless $j\in [a-(r_\sigma - r_\eta), b]$.

Since
\begin{equation}
\i_\sigma^\circ(\ic_\p(\Delta))^i =
\bigoplus_{\eta\in F(\sigma)\setminus\{\sigma\}}
\ic_\p(\Delta)(\eta)_{A(\sigma)}^i\;,
\end{equation}
(2) is a consequence of (1) for
$\eta\in F(\sigma)\setminus\{\sigma\}$ which are true by
the assumption of the induction.
Since
\begin{equation}
\ic_\p(\Delta)(\sigma)^i =
\grt^{\geq\p(\sigma)+1}(\i_\sigma^\circ\ic_\p(\Delta)[-1])^i\;,
\end{equation}
(1) follows from (2).
Since
\begin{equation}
\i_\sigma^*(\ic_\p(\Delta))^i =
\i_\sigma^\circ(\ic_\p(\Delta))^i\oplus\ic_\p(\Delta)(\sigma)^i
\end{equation}
for every $i\in\Z$, (3) follows from (1) and (2).

Since
\begin{equation}
\Gamma(\ic_\p(\Delta))^i =
\bigoplus_{\sigma\in\Delta}\ic_\p(\Delta)(\sigma)_A^i
\end{equation}
for every $i$, the last assertion is a consequence of
(1) for all $\sigma\in\Delta\setminus\{\Zero\}$.
\QED

%
%

\begin{Cor}
\label{Cor 2.12}
For any perversity $\p$ on $\Delta$, $\D(\ic_\p(\Delta))^\bullet$ is
quasi-isomoprhic to $\ic_{-\p}(\Delta)^\bullet$.
\end{Cor}

\Proof
It is sufficient to show that $\D(\ic_\p(\Delta))^\bullet$
satisfies the conditions of the theorem for the perversity $-\p$.
(1) is satisfied since $\D(\ic_\p(\Delta))(\Zero)^0 =\Q$ and
$\D(\ic_\p(\Delta))(\Zero)^i =\{0\}$ for $i\not= 0$ by the
definition of $\D$.

We check the conditions (2) and (3) for each $\sigma$ in
$\Delta\setminus\{\Zero\}$.

Since
$\D(\ic_\p(\Delta))(\sigma)^\bullet =
\det(\sigma)\otimes\d_\sigma(\i_\sigma^*\ic_\p(\Delta))[-r_\sigma]^\bullet$,
the condition (3) of the theorem and Lemma~\ref{Lem 1.6} imply
that $\grt_{\leq-\p(\sigma)}(\D(\ic_\p(\Delta))(\sigma))^\bullet$
is acyclic, i.e., (2) for $-\p$.

By Lemma~\ref{Lem 2.2}, $\ic_\p(\Delta)(\sigma)^\bullet$ is
quasi-isomorphic to $\D(\D(\ic_\p(\Delta)))(\sigma)^\bullet$.
Hence $\grt_{\leq\p(\sigma)}(\D(\D(\ic_\p(\Delta)))(\sigma))^\bullet$
is acyclic by (2) for $\ic_\p(\Delta)(\sigma)^\bullet$.
Since
\begin{equation}
\det(\sigma)\otimes\d_\sigma(\i_\sigma^*\D(\ic_\p(\Delta)))^\bullet =
\D(\D(\ic_\p(\Delta)))(\sigma)[r_\sigma]^\bullet\;,
\end{equation}
Lemma~\ref{Lem 1.6} implies that
$\grt^{\geq-\p(\sigma)}(\i_\sigma^*\D(\ic_\p(\Delta)))^\bullet$
is acyclic, i.e., (3) for $-\p$.
\QED

Two homomorphisms $f, g\cl L^\bullet\ra K^\bullet$ in
$\CGEM(\Delta)$ are said to be homotopic if there exists a
collection of homomorphisms
$\{u^i\cl L^i\ra K^{i-1}\scl i\in\Z\}$ in $\GEM(\Delta)$
such that
\begin{equation}
f^i - g^i = d_K^{i-1}\cdot u^i + u^{i+1}\cdot d_L^i
\end{equation}
for every $i\in\Z$.
If $f$ and $g$ are homotopic, then
$f(\sigma),g(\sigma)\cl L(\sigma)^\bullet\ra K(\sigma)^\bullet$
and
$\i_\sigma^*(f),\i_\sigma^*(g)\cl\i_\sigma^*L^\bullet\ra\i_\sigma^*K^\bullet$
for $\sigma\in\Delta$ as well as
$\Gamma(f),\Gamma(\gamma)\cl\Gamma(L)^\bullet\ra\Gamma(K)^\bullet$
are homotopic as complexes in abelian categories.
Actually, it is sufficient to take
$\{u^i(\sigma/\sigma)\}$, $\{\i_\sigma^*(u^i)\}$ and $\{\Gamma(u^i)\}$,
respectively.
In particular, the maps of the cohomologies induced by the
two homomorphisms of the complexes are respectively equal.

We give here some elementary lemmas on the quasi-isomorphism
property in $\CGEM(\Delta)$.

%
%

\begin{Lem}
\label{Lem 2.13}
Let $f_1\cl L_1^\bullet\ra L_2^\bullet$ be a quasi-isomorphism
and $f_2\cl L_1^\bullet\ra L_3^\bullet$ a homomorphism
in $\CGEM(\Delta)$.
Then there exist $L_4^\bullet$ in $\CGEM(\Delta)$,
a quasi-isomorphism $g_1\cl L_3^\bullet\ra L_4^\bullet$ and
a homomorphism $g_2\cl L_2^\bullet\ra L_4^\bullet$ such that
the homomorphisms $g_2\cdot f_1$ and $g_1\cdot f_2$ are homotopic.
If $f_2$ is a quasi-isomorphism, then so is $g_2$.
\end{Lem}

%
%

\begin{Lem}
\label{Lem 2.14}
Let $g_1\cl L_3^\bullet\ra L_4^\bullet$ be a quasi-isomorphism
and $g_2\cl L_2^\bullet\ra L_4^\bullet$ a homomorphism
in $\CGEM(\Delta)$.
Then there exist $L_1^\bullet$ in $\CGEM(\Delta)$,
a quasi-isomorphism $f_1\cl L_1^\bullet\ra L_2^\bullet$ and
a homomorphism $f_2\cl L_1^\bullet\ra L_3^\bullet$ such that
the homomorphisms $g_2\cdot f_1$ and $g_1\cdot f_2$ are homotopic.
If $g_2$ is a quasi-isomorphism, then so is $f_2$.
\end{Lem}

We get these lemmas by setting $L_4^\bullet$ and
$L_1^\bullet$ the mapping cones of the homomorphisms
$L_1^\bullet\ra L_2^\bullet\oplus L_3^\bullet$ and
$L_2^\bullet\oplus L_3^\bullet\ra L_4^\bullet$, respectively.

By applying these lemmas, we get the following equivalent
conditions.

%
%

\begin{Lem}
\label{Lem 2.15}
For $L^\bullet, K^\bullet$ in $\CGEM(\Delta)$, the following
conditions are equivalent.

(1) $L^\bullet$ is quasi-isomorphic to $K^\bullet$.

(2) There exists $J^\bullet$ in $\CGEM(\Delta)$ and
quasi-isomorphisms $L^\bullet\ra J^\bullet$ and
$K^\bullet\ra J^\bullet$.

(3) There exists $I^\bullet$ in $\CGEM(\Delta)$ and
quasi-isomorphisms $I^\bullet\ra L^\bullet$ and
$I^\bullet\ra K^\bullet$.
\end{Lem}

%
%

\begin{Lem}
\label{Lem 2.16}
Let $f_1\cl L_1^\bullet\ra K_1^\bullet$ be a homomorphism
in $\CGEM(\Delta)$.
Assume that $L_2^\bullet, K_2^\bullet$ in $\CGEM(\Delta)$ are
quasi-isomorphic to $L_1^\bullet$ and $K_1^\bullet$, respectively.
Then there exist $K_3^\bullet$ in $\CGEM(\Delta)$, a quasi-isomorphism
$K_2^\bullet\ra K_3^\bullet$ and a homomorphism
$f_2\cl L_2^\bullet\ra K_3^\bullet$ such that the diagrams
\begin{equation}
\begin{array}{ccc}
   \H^i(\Gamma(L_1)^\bullet)
 & \lra & H^i(\Gamma(K_1)^\bullet) \\
    \wr |  &  & \wr | \\
   \H^i(\Gamma(L_2)^\bullet)
 & \lra & H^i(\Gamma(K_3)^\bullet) \\
\end{array}
\end{equation}
of the cohomologies are commutative for all $i\in\Z$.

There exist also $L_3^\bullet$ in $\CGEM(\Delta)$, a quasi-isomorphism
$L_3^\bullet\ra L_2^\bullet$ and a homomorphism
$f_3\cl L_3^\bullet\ra K_2^\bullet$ with the similar compatibility
with $f_1$.
\end{Lem}

\Proof
By Lemma~\ref{Lem 2.15}, there exists $L_0^\bullet$ and
quasi-isomorphisms $g_1\cl L_0^\bullet\ra L_1^\bullet$ and
$g_2\cl L_0^\bullet\ra L_2$.
By applying Lemma~\ref{Lem 2.14} for $f_1\cdot g_1$ and
$g_2$, we get $K_0^\bullet$ with a quasi-isomorphism
$K_1^\bullet\ra K_0^\bullet$ and a homomorphism
$h_0\cl L_2^\bullet\ra K_0^\bullet$ with the compatibility
condition.
Since $K_0^\bullet$ is quasi-isomorphic to $K_2$,
there exists $K_3^\bullet$ and quasi-isomorphisms
$h_1\cl K_0^\bullet\ra K_3^\bullet$ and
$K_2^\bullet\ra K_3^\bullet$ by Lemma~\ref{Lem 2.15}.
It is sufficient to set $f_3 := h_1\cdot h_0$.

The second assertion is proved similarly.
\QED


\section{The algebraic theory on toric varieties}



\setcounter{equation}{0}

In this section, we construct a functor from the category of
graded exterior modules to that of complexes on the toric
variety associated to the fan.
For a finite fan $\Delta$, we denote by $Z(\Delta)$ the associated
toric variety defined over $\Q$ (cf.\cite{Oda1}).

We start with the case of an affine toric variety.
Assume that $\Delta$ is $F(\pi)$, i.e. the set of all faces of
a cone $\pi$ in $N_\R$.

We denote by $\Q[M]$ the group ring $\bigoplus_{m\in M}\Q\e(m)$
defined by $\e(m)\e(m') =\e(m + m')$ for $m, m'\in M$ and $\e(0) = 1$.
This $\Q$-algebra has a grading in the free $\Z$-module $M$.
For a subset $U$ of $M$, we denote
$\Q[U] :=\bigoplus_{m\in U}\Q\e(m)$.
Note that $1\not\in\Q[U]$ if $0\not\in U$.

For the subsemigroup $M\cap\pi^\vee\subset M$,
we denote by $S(\pi)$ the $M$-graded $\Q$-subalgebra
$\Q[M\cap\pi^\vee]$ of $\Q[M]$.
Then the affine toric variety $Z(F(\pi))$ is equal to
$\Spec S(\pi)$.
The algebraic torus $T_N$ is equal to $\Spec\Q[M]$ and the reduced
complement $Z(F(\pi))\setminus T_N$ is defined by the ideal
$J(\pi) :=\Q[M\cap\inte\pi^\vee]$.

The {\em logarithmic de Rham complex}
$\Omega_{S(\pi)}(\log J(\pi))^\bullet$ is defined as follows.

We set
\begin{equation}
\Omega_{S(\pi)}(\log J(\pi))^1 := S(\pi)\otimes M
\end{equation}
and
\begin{equation}
\Omega_{S(\pi)}(\log J(\pi))^i :=
\bigwedge^i\Omega_{S(\pi)}(\log J(\pi))^1 =
S(\pi)\otimes\bigwedge^iM
\end{equation}
for $0\leq i\leq r$, where the exterior powers are taken as
an $S(\pi)$-module and as a $\Z$-module, respectively.
These are clearly free $S(\pi)$-modules.
By the notation $A^* = A(M_\Q) =\bigwedge^\bullet M_\Q$,
the direct sum $\bigoplus_{i=0}^r\Omega_{S(\pi)}(\log J(\pi))^i$
is equal to the $M$-graded free $S(\pi)$-module
\begin{equation}
S(\pi)\otimes_\Q A^*
=\bigoplus_{m\in M\cap\pi^\vee}\Q\e(m)\otimes_\Q A^*\;.
\end{equation}
The $\Q$-endomorphism $\partial$ of this $S(\pi)$-module
is defined to be the $M$-homogeneous morphism such
that the the restriction to the component $\Q\e(m)\otimes_\Q A^*$
is $1\otimes d_m$ for each $m\in M\cap\pi^\vee$, where $d_m$
is the left multiplication of $m$.
Since
$\partial(\Omega_{S(\pi)}(\log J(\pi))^i)\subset\Omega_{S(\pi)}(\log
J(\pi))^{i+1}$
for each $i$, $\Omega_{S(\pi)}(\log J(\pi))^\bullet$ is a
$\partial$-complex of $M$-graded $\Q$-vector spaces.

For each face $\sigma$ of $\pi$, $\pi^\vee\cap\sigma^\bot$ is a
face of the dual cone $\pi^\vee\subset M_\R$.
Furthermore, it is known that the correspondence
$\sigma\mapsto\pi^\vee\cap\sigma^\bot$ defines a bijection from
$F(\pi)$ to $F(\pi^\vee)$ \cite[Prop.A.6]{Oda1}.

For each $\sigma\in F(\pi)$, we denote by $P(\pi;\sigma)$ the
$M$-homogeneous prime ideal
\begin{equation}
\Q[M\cap(\pi^\vee\setminus\sigma^\bot)] =
\bigoplus_{m\in M\cap(\pi^\vee\setminus\sigma^\bot)}\Q\e(m)
\end{equation}
of $S(\pi)$.
The $M$-homogeneous quotient ring $S(\pi)/P(\pi;\sigma)$ is denoted
by $S(\pi;\sigma)$. We denote the image of $\e(m)$ in $S(\pi;\sigma)$
for $m\in M\cap\pi^\vee\cap\sigma^\bot$ also by $\e(m)$.
Since $M[\sigma] = M\cap\sigma^\bot$, we have a description
$S(\pi;\sigma) =\bigoplus_{m\in M[\sigma]\cap\pi^\vee}\Q\e(m)$.

Let $J(\pi;\sigma)$ be the ideal
$\Q[M[\sigma]\cap\relint(\pi^\vee\cap\sigma^\bot)]$ of $S(\pi;\sigma)$.
The $\partial$-complex $\Omega_{S(\pi;\sigma)}(\log J(\pi;\sigma))^\bullet$
is defined to be $S(\pi;\sigma)\otimes_\Q A^*[\sigma]$ with the
$M$-homogeneous $\Q$-homomorphism $\partial$ defined similarly as above,
where $A^*[\sigma] =\bigwedge^\bullet M[\sigma]_\Q$.
Note that $\Omega_{S(\pi;\Zero)}(\log J(\pi;\Zero))^\bullet$ is equal to
$\Omega_{S(\pi)}(\log J(\pi))^\bullet$.

We denote by $\Coh(S(\pi))$ the category:
\begin{description}
\item[object:]  A finitely generated $M$-graded $S(\pi)$-module.
\item[morphism:] An $M$-homogeneous $S(\pi)$-homomorphism
of $M$-degree zero.
\end{description}

A $\Q$-homomorphism $\delta\cl F\ra G$ of $S(\pi)$-modules is said to be
a {\em differential operator of order one} if the map
$(\delta\cdot f - f\cdot\delta)\cl F\ra G$ defined by
$(\delta\cdot f - f\cdot\delta)(x) :=\delta(fx) - f(\delta(x))$
is an $S(\pi)$-homomorphism for every $f\in S(\pi)$.
We denote by $\CCohDiff(S(\pi))$ the category:
\begin{description}
\item[object:] A finite $\partial$-complex $F^\bullet$ such that
$F^i$'s are in $\Coh(S(\pi))$ and $\partial$ is $M$-homogeneous
of $M$-degree zero and is a differential operator of order one.
\item[morphism:] An $M$-homogeneous $S(\pi)$-homomorphism
of $M$-degree zero.
\end{description}

We construct a functor $\Lambda_{S(\pi)}$ from the category
$\GEM(F(\pi))$ of graded exterior modules on $F(\pi)$ to
this category $\CCohDiff(S(\pi))$.

Let $\rho$ be a cone in $F(\pi)$.
Recall that an object $V$ of $\GM(A(\rho))$ is a
finitely generated graded $A(\rho)$-module and the $A$-module
$V_A$ has a structure of a free $A^*[\rho]$-module
(cf. Lemma~\ref{Lem 1.3}).
Each $m\in M[\rho] = M\cap\rho^\bot$ is a homogeneous element of
$A^*[\rho]$ of degree one.
We denote by $d_m$ the left operation of $m$ on $V_A$.
Then $d_m^2 = 0$ since $m\wedge m = 0$.
For each $m\in M[\rho]$, we denote by $V_A(m)^\bullet$
the $\partial$-complex defined by $V_A(m)^i :=(V_A)_i$
for each $i\in\Z$ and $\partial := d_m$.

We set $\Lambda_{S(\pi)}^\rho(V)^i := S(\pi;\rho)\otimes_\Q(V_A)_i$
for each integer $i$.
We define the $\partial$-complex $\Lambda_{S(\pi)}^\rho(V)^\bullet$
by
\begin{equation}
\Lambda_{S(\pi)}^\rho(V)^\bullet :=S(\pi;\rho)\otimes_\Q V_A =
\bigoplus_{m\in M[\rho]\cap\pi^\vee}\Q\e(m)\otimes_\Q V_A(m)^\bullet\;,
\end{equation}
i.e., the its $m$-component of $\partial$ is $1_{\Q\e(m)}\otimes d_m$
for every $m\in M[\rho]\cap\pi^\vee$.

In order to check that $\partial$ is a differential operator of
order one, it is sufficient to show the $S(\pi;\rho)$-linearlity of
$\partial\cdot \e(m_0) -\e(m_0)\cdot \partial$ for $m_0\in M\cap\pi^\vee$.
Since $P(\pi;\rho)$ is the annihilator of $S(\pi;\rho)$,
$\partial\cdot \e(m_0) -\e(m_0)\cdot \partial = 0$ if $m_0\not\in M[\rho]$.
Assume $m_0, m_1\in M[\rho]\cap\pi^\vee$ and $x\in V_A$.
For any $m\in M[\rho]\cap\pi^\vee$, we have
\begin{eqnarray*}
 &   & (\partial\cdot \e(m_0) -\e(m_0)\cdot \partial)(\e(m)\e(m_1)\otimes x) \\
 & = & \partial(\e(m + m_0 + m_1)\otimes x) -\e(m_0)\partial(\e(m + m_1)\otimes
x) \\
 & = & \e(m + m_0 + m_1)\otimes(m + m_0 + m_1)x
       -\e(m + m_0 + m_1)\otimes(m + m_1)x \\
 & = & \e(m + m_0 + m_1)\otimes m_0x \\
 & = & \e(m + m_0 + m_1)\otimes(m_0 + m_1)x
       -\e(m + m_0 + m_1)\otimes m_1x \\
 & = & \e(m)(\partial\cdot \e(m_0) -\e(m_0)\cdot \partial)(\e(m_1)\otimes x)\;.
\end{eqnarray*}
Hence $\partial\cdot \e(m_0) -\e(m_0)\cdot \partial$ is an
$S(\pi;\rho)$-homomorphism.

%
%

\begin{Prop}
\label{Prop 3.1}
Let $\rho$ be in $F(\pi)$ and $V$ in $\GM(A(\rho))$.
Then $\Lambda_{S(\pi)}^\rho(V)^\bullet$ is isomorphic to a finite
direct sum of subcomplexes which are isomorphic to dimension shifts
of $\Omega_{S(\pi;\rho)}(\log J(\rho))^\bullet$.
\end{Prop}

\Proof
By Lemma~\ref{Lem 1.3}, $V_A$ is a free $A^*[\rho]$-module.
Let $\{x_1,\cdots, x_s\}$ be a homogeneous basis.
By the definition of $\partial$, the decomposition
\begin{equation}
\Lambda_{S(\pi)}^\rho(V)^\bullet =
\bigoplus_{i=1}^sS(\pi;\rho)\otimes_\Q A^*[\rho]x_i
\end{equation}
is a direct sum of subcomplexes.
Furthermore, there exists an isomorphism
\begin{equation}
\Omega_{S(\pi;\rho)}(\log J(\rho))[-\deg x_i]^\bullet\lra
(S(\pi;\rho)\otimes_\Q A^*[\rho]x_i)^\bullet
\end{equation}
for each $i$.
\QED

For a homomorphism $f\cl V\ra W$ in $\GM(A(\rho))$, the
$S(\pi)$-homomorphism
$\Lambda_{S(\pi)}^\rho(f)\cl\Lambda_{S(\pi)}^\rho(V)^\bullet\ra%
\Lambda_{S(\pi)}^\rho(W)^\bullet$ of $\partial$-complexes is defined
by $\Lambda_{S(\pi)}^\rho(f) = 1_{S(\pi)}\otimes f_A$.
Since $f_A$ is an $A^*[\rho]$-homomorphism, $\Lambda_{S(\pi)}^\rho(f)$
commutes with the coboudary maps of $\Lambda_{S(\pi)}^\rho(V)^\bullet$
and $\Lambda_{S(\pi)}^\rho(W)^\bullet$.

%
%

\begin{Prop}
\label{Prop 3.2}
Let $\sigma,\rho$ be cones in $F(\pi)$ with $\sigma\prec\rho$.
For $V$ in $\GM(A(\sigma))$, there exists a natural isomorphism
\begin{equation}
\Lambda_{S(\pi)}^\sigma(V)^\bullet\otimes_{S(\pi;\sigma)}S(\pi;\rho)\simeq
\Lambda_{S(\pi)}^\rho(V_{A(\rho)})^\bullet
\end{equation}
of $M$-graded $\partial$-complexes.
\end{Prop}

\Proof
$\Lambda_{S(\pi)}^\sigma(V)^\bullet\otimes_{S(\pi;\sigma)}S(\pi;\rho)$
is actually an $M$-graded $\partial$-complex since it is the quotient
of $\Lambda_{S(\pi)}^\sigma(V)^\bullet$
by the $M$-homogeneous subcomplex
$P(\pi;\rho)\Lambda_{S(\pi)}^\sigma(V)^\bullet$.
Both sides are equal to $S(\pi;\rho)\otimes_\Q V_A$ by the
identification $(V_{A(\rho)})_A = V_A$.
The differential operators commute with the identification
since the $m$-components of them for $m\in M[\rho]\cap\pi^\vee$
are both the operation $d_m$.
\QED

For $L$ in $\GEM(F(\pi))$, we define the $\partial$-complex
$\Lambda_{S(\pi)}(L)^\bullet$ in $\CCohDiff(S(\pi))$ by
\begin{equation}
\Lambda_{S(\pi)}(L)^\bullet :=
\bigoplus_{\sigma\in F(\pi)}\Lambda_{S(\pi)}^\rho(L(\rho))^\bullet\;,
\end{equation}
where the coboundary map $\partial$ is also defined as the direct sum.

Let $f\cl L\ra K$ be a morphism in $\GEM(F(\pi))$.
We define an $M$-homogeneous $S(\pi)$-homomorphism
$\Lambda_{S(\pi)}(f)\cl\Lambda_{S(\pi)}(L)^\bullet\ra%
\Lambda_{S(\pi)}(K)^\bullet$ of $\partial$-complexes as follows.

For $\sigma,\rho\in F(\pi)$, the $(\sigma,\rho)$-component of
the morphism
\begin{equation}
\begin{array}{ccccc}
\Lambda_{S(\pi)}(f) & \cl & \Lambda_{S(\pi)}(L)^\bullet
 & \lra  & \Lambda_{S(\pi)}(K)^\bullet \\
                    &     & \| &    & \| \\
                    &     &
\bigoplus_{\sigma\in F(\pi)}\Lambda_{S(\pi)}^\sigma(L(\sigma))^\bullet
 &  &
\bigoplus_{\rho\in F(\pi)}\Lambda_{S(\pi)}^\rho(K(\rho))^\bullet
\end{array}
\end{equation}
is defined to be the composite of the natural surjection
\begin{equation}
\Lambda_{S(\pi)}^\sigma(L(\sigma))^\bullet\ra
\Lambda_{S(\pi)}^\sigma(L(\sigma))^\bullet%
\otimes_{S(\pi;\sigma)}S(\pi;\rho)
= \Lambda_{S(\pi)}^\rho(L(\sigma)_{A(\rho)})^\bullet
\end{equation}
and
\begin{equation}
1_{S(\pi;\rho)}\otimes f(\sigma/\rho)_A\cl
\Lambda_{S(\pi)}^\rho(L(\sigma)_{A(\rho)})^\bullet\lra
\Lambda_{S(\pi)}^\rho(K(\rho))^\bullet
\end{equation}
if $\sigma\prec\rho$ and is defined to be the zero map otherwise.
Then $\Lambda_{S(\pi)}$ is a covariant functor from $\GEM(\Delta)$
to $\CCohDiff(S(\pi))$.

Let $k$ be a field.
For a topological space $X$, we denote by $k_X$ the sheaf of
rings with the constant stalk $k$.
A $k_X$-module on $X$ is called a $k$-sheaf and
a $k_X$-homomorphism of $k_X$-modules is called simply a
$k$-homomorphism.
In this paper, we treat only the cases $k =\Q$ and $k =\C$.

Let $\calO_X$ be a sheaf of commutative $k$-algebras.
A $k$-homomoprhism $f\cl F\ra G$ of $\calO_X$-modules $F, G$
is said to be a differential operator of order one if the
$k$-homomorphism $f(U)\cl F(U)\ra G(U)$ is a differential
operator of order one of $\calO_X(U)$-modules
for every open subset $U\subset X$.

For a finitely generated $S(\pi)$-module $E$, we denote by $\calF(E)$
the associated coherent sheaf on the affine scheme $\Spec S(\pi)$.
If $d\cl E\ra G$ is a differential operator of order one of
finitely generated $S(\pi)$-modules, then it defines a differential
operator $\calF(d)\cl\calF(E)\ra\calF(G)$ of order one of
$\calO_{Z(F(\pi))}$-modules on the affine toric variety $Z(F(\pi))$.
If $(E^\bullet, \partial = (\partial_E^i))$ is a $\partial$-complex
such that each $E^i$ is a finitely generated $S(\pi)$-module and
$\partial$ is a differential operator of order one, then we
denote by $\calF(E)^\bullet$ the $\partial$-complex on
$Z(F(\pi))$ with the coboudary map
$\partial = (\calF(\partial_E^i))$.

Let $\mu$ be an element of $F(\pi)$.
Take an element $m\in M\cap\relint(\pi^\vee\cap\mu^\bot)$.
Since $\mu^\vee =\pi^\vee +\R(-m)$ \cite[Cor.A.7]{Oda1},
$S(\mu) :=\Q[M\cap\mu^\vee]$ is equal to the localization
$S(\pi)[\e(m)^{-1}]$.
For any $\rho\in F(\pi)$, we see easily that
$S(\pi;\rho)\otimes_{S(\pi)}S(\mu)$ is equal to $S(\mu;\rho)$
if $\rho\prec\mu$ and is $\{0\}$ otherwise.

%
%

\begin{Lem}
\label{Lem 3.3}
Let $\mu$ be an element of $F(\pi)$.
For $\rho\in F(\mu)$ and $V$ in $\GM(A(\rho))$, the localization
$\Lambda_{S(\pi)}^\rho(V)^\bullet\otimes_{S(\pi)}S(\mu)$
is equal to $\Lambda_{S(\mu)}^\rho(V)^\bullet$.

For $L$ in $\GEM(F(\pi))$, the localization
$\Lambda_{S(\pi)}(L)^\bullet\otimes_{S(\pi)}S(\mu)$
of $\partial$-complex is equal to
$\Lambda_{S(\mu)}(L|{F(\mu)})^\bullet$, where $L|{F(\mu)}$ is the
restriction of $L$ to $F(\mu)$.
\end{Lem}

\Proof
The first equality is clear as $S(\mu)$-modules.
The coboundary maps $\partial$ are compatible with the
inclusion
$\Lambda_{S(\pi)}^\rho(V)^\bullet\subset\Lambda_{S(\mu)}^\rho(V)^\bullet$
since they are defined as the direct sums of
$\Q\e(m)\otimes_\Q V_A(m)^\bullet$'s.
Since differential operators are extended to localizations uniquely,
they are equal as $\partial$-complexes.

Since $\Lambda_{S(\pi)}(L)^\bullet$ and $\Lambda_{S(\mu)}(L)^\bullet$
are defined as the direct sums of
$\Lambda_{S(\pi)}^\rho(L(\rho))^\bullet$ for $\rho\in F(\pi)$ and
$\Lambda_{S(\mu)}^\rho(L(\rho))^\bullet$ for $\rho\in F(\mu)$,
respectively, the second assertion follows from the first.
\QED

Let $\Delta$ be a finite fan of $N_\R$.
The toric variety $Z(\Delta)$ has the affine open covering
$\bigcup_{\pi\in\Delta}Z(F(\pi))$.
The reduced divisor $D(\Delta) := Z(\Delta)\setminus T_N$
is defined by the ideal $J(\pi)\subset S(\pi)$ on each affine
open subscheme $Z(F(\pi))$.

The logarithmic de Rham complex
$\Omega_{Z(\Delta)}(\log D(\Delta))^\bullet$ is defined by
\begin{equation}
\Omega_{Z(\Delta)}(\log D(\Delta))^\bullet :=
\calO_{Z(\Delta)}\otimes_\Q A^*\;.
\end{equation}
We define the coboudary map $\partial = (\partial^i\cl i\in\Z)$
\begin{equation}
\partial^i\cl\Omega_{Z(\Delta)}(\log D(\Delta))^i\lra
\Omega_{Z(\Delta)}(\log D(\Delta))^{i+1}
\end{equation}
so that the restriction to $Z(F(\pi))$ is equal to $\partial$
of $\calF(\Omega_{S(\pi)}(\log J(\pi)))^\bullet$ for each
$\pi\in\Delta$.

Although it is common to write a logarithmic de Rham complex
as $\Omega_X^\bullet(\log D)$, we put the dot at the right end
as $\Omega_X(\log D)^\bullet$ for the compatibility with the other
notation in this paper.

For each $\sigma\in F(\rho)$, the subscheme $\Spec S(\rho;\sigma)$
of $Z(F(\rho))$ is denoted by $X(\rho;\sigma)$.
In particular, $X(\rho;\rho)$ is the algebraic torus
$T_{N[\rho]} :=\Spec\Q[M[\rho]]$ of dimension $r - r_\rho$.

In order to simplify the notation, we set $T := T_N$ and
$T[\rho] := T_{N[\rho]}$ for each $\rho\in\Delta$.

Then the toric variety $Z(\Delta)$ is decomposed as the
disjoint union
\begin{equation}
\bigcup_{\rho\in\Delta}T[\rho]
\end{equation}
of $T$-orbits \cite[Prop.1.6]{Oda1}.

For each $\sigma\in\Delta$, we denote by $X(\Delta;\sigma)$ or
simply $X(\sigma)$ the union of $X(\rho;\sigma)$ for $\rho\in\Delta$
with $\sigma\prec\rho$. $X(\sigma)$ is a $T$-invariant irreducible
closed subvariety of $Z(\Delta)$.

For each $\sigma\in\Delta$, let $N[\sigma] := N/N(\sigma)$.
For each $\rho\in\Delta$ with $\sigma\prec\rho$, we denote by
$\rho[\sigma]$ the image of $\rho$ in $N[\sigma]_\R = N_\R/N(\sigma)_\R$.
Then $\Delta[\sigma] :=\{\rho[\sigma]\cl\rho\in\Delta,\sigma\prec\rho\}$
is a fan of $N[\sigma]_\R$.
It is known that $X(\sigma)\subset Z(\Delta)$ is equal to
the toric variety $Z(\Delta[\sigma])$ with the torus $T[\sigma]$
\cite[Cor.1.7]{Oda1}.
The reduced complement $X(\sigma)\setminus T[\sigma]$ is denoted by
$D(\Delta;\sigma)$ or $D(\sigma)$.

We define
\begin{equation}
\Omega_{X(\sigma)}(\log D(\sigma))^\bullet :=
\calO_{X(\sigma)}\otimes_\Q A^*[\sigma]\;
\end{equation}
and $\partial$ of it is defined so that the restriction to the
affine open subscheme $Z(F(\rho))$ is equal to that of
$\calF(\Omega_{S(\rho;\sigma)}(\log J(\rho;\sigma)))^\bullet$ for
every $\rho\in\Delta$ with $\sigma\prec\rho$.

For each $V\in\GM(A(\sigma))$, we define the
$\partial$-complex
$\Lambda_{Z(\Delta)}^\sigma(V)^\bullet$ on $Z(\Delta)$ by
\begin{equation}
\Lambda_{Z(\Delta)}^\sigma(V)^\bullet :=
\calO_{X(\sigma)}\otimes_\Q V_A\;.
\end{equation}
The coboudary map $\partial$ is defined so that the restriction to
each open set $Z(F(\rho))$ is equal to that of
$\calF(\Lambda_{S(\rho)}^\sigma(V))^\bullet$ for every
$\rho\in\Delta$ with $\sigma\prec\rho$.

Note that $M$-gradings of $\Lambda_{S(\rho)}^\sigma(V)^\bullet$
for $\rho\in\Delta$ induce a natural $T$-action on the $\partial$-complex
$\Lambda_{Z(\Delta)}^\sigma(V)^\bullet$.

The following propositions follow from Propositions~\ref{Prop 3.1}
and \ref{Prop 3.2}, respectively.

%
%

\begin{Prop}
\label{Prop 3.4}
Let $\rho$ be in $\Delta$ and $V$ in $\GM(A(\rho))$.
Then $\Lambda_{Z(\Delta)}^\rho(V)^\bullet$ is isomorphic to a
finite direct sum of subcomplexes which are isomorphic to dimension
shifts of $\Omega_{X(\rho)}(\log D(\rho))^\bullet$.
\end{Prop}

%
%

\begin{Prop}
\label{Prop 3.5}
Let $\sigma,\rho$ be cones in $\Delta$ with $\sigma\prec\rho$.
For $V$ in $\GM(A(\sigma))$, there exists a natural
$T$-equivariant isomorphism
\begin{equation}
\Lambda_{Z(\Delta)}^\sigma(V)^\bullet\otimes_{\calO_{X(\sigma)}}
\calO_{X(\rho)}\simeq
\Lambda_{Z(\Delta)}^\rho(V_{A(\rho)})^\bullet
\end{equation}
of $\partial$-complexes.
\end{Prop}

Let $V^\bullet$ be in $\CGM(A(\sigma))$.
Then the bicomplex
$\Lambda_{Z(\Delta)}^\sigma(V)^{\bullet,\bullet}$ is defined by
\begin{equation}
\Lambda_{Z(\Delta)}^\sigma(V)^{i,j} :=
\Lambda_{Z(\Delta)}^\sigma(V^i)^j
\end{equation}
for $i, j\in\Z$ and $d_1 :=d$ and $d_2 :=\partial$.
We denote by $\Lambda_{Z(\Delta)}^\sigma(V)^\bullet$ the
associated single complex and by $\delta$ the
coboundary map.

For each object $L$ of $\GEM(\Delta)$, we set
\begin{equation}
\Lambda_{Z(\Delta)}(L)^\bullet :=
\bigoplus_{\sigma\in\Delta}%
\Lambda_{Z(\Delta)}^\sigma(L(\sigma))^\bullet\;.
\end{equation}
For a morphism $f\cl L\ra K$ in $\GEM(\Delta)$, the
$T$-equivariant homomorphism
\begin{equation}
\Lambda_{Z(\Delta)}(f)\cl\Lambda_{Z(\Delta)}(L)^\bullet%
\lra\Lambda_{Z(\Delta)}(K)^\bullet
\end{equation}
is defined naturally.
Then $\Lambda_{Z(\Delta)}$ is a covariant functor from
$\GEM(\Delta)$ to the category $\CCohDiff(Z(\Delta))$ which
is defined naturally as the globalization of
$\CCohDiff(S(\pi))$.

Let $L^\bullet$ be a $d$-complex in $\CGEM(\Delta)$.
Then the bicomplex $\Lambda_{Z(\Delta)}(L)^{\bullet,\bullet}$
is defined by
\begin{equation}
\Lambda_{Z(\Delta)}(L)^{i,j} :=\Lambda_{Z(\Delta)}(L^i)^j
\end{equation}
for $i,j\in\Z$.
Note that $d_1 := d$ of this bicomplex is a $\calO_{Z(\Delta)}$-homomorphism
and $d_2 :=\partial$ is a differential operator of order one.
If there is no danger of confusion, we denote by
$\Lambda_{Z(\Delta)}(L)^\bullet$ the associated single complex.
The coboudary map, which we denote by $\delta$, is a differential
operator of order one.
For each integer $j$, we denote by $\Lambda_{Z(\Delta)}(L)_j^\bullet$
the $d$-complex $\Lambda_{Z(\Delta)}(L)^{\bullet,j}$.
Then $\Lambda_{Z(\Delta)}(L)_j^\bullet$ is a finite
$d$-complex in the category of coherent $\calO_Z$-modules.

If the fan $\Delta$ is complete, then the toric variety $Z(\Delta)$
is complete.
For the functor $\Gamma\cl\GEM(\Delta)\ra\GM(A)$ defined in
Section~1, we get the following lemma.

%
%

\begin{Lem}
\label{Lem 3.6}
Assume that $\Delta$ is a complete fan.
Let $L^\bullet$ be an object of $\CGEM(\Delta)$.
For any integers $p, q$, we have an isomorphism
\begin{equation}
\H^p(Z(\Delta),\Lambda_Z(L)_q^\bullet)\simeq
\H^p(\Gamma(L)^\bullet)_q
\end{equation}
of finite dimensional $\Q$-vector spaces, where
the lefthand side is the hepercohomology group of the
complex of coherent sheaves on $Z(\Delta)$.
\end{Lem}

\Proof
For each $\sigma\in\Delta$, we have
$\H^0(X(\sigma),\calO_{X(\sigma)}) =\Q$ and
$\H^p(X(\sigma),\calO_{X(\sigma)}) = \{0\}$ for $p > 0$, since
$X(\sigma)$ is a complete toric variety \cite[Cor.2.8]{Oda1}.
Since $\Lambda_Z(L)_q^p =\Lambda_Z(L)^{p,q}$ is a direct sum for
$\sigma\in\Delta$ of free $\calO_{X(\sigma)}$-modules for every
$p\in\Z$, the hypercohomology $\H^p(Z(\Delta),\Lambda_Z(L)_q^\bullet)$
is equal to the $p$-th cohomology of the complex
$\Gamma(Z,\Lambda_Z(L)_q)^\bullet$ of $\Q$-vector spaces.
For any $\sigma\in\Delta$ and $p, q\in\Z$, we have
\begin{equation}
\Gamma(Z,\Lambda_Z^\sigma(L(\sigma))_q)^p
= (L(\sigma)_A^p)_q\;.
\end{equation}
Hence $\Gamma(Z,\Lambda_Z(L)_q)^\bullet$ is isomorphic
to $\Gamma(L)_q^\bullet$ as a complex of $\Q$-vector spaces.
\QED


\section{The analytic theory on toric varieties}



\setcounter{equation}{0}

For any $\Q$-algebra $B$ and any $\Q$-scheme $X$, we denote by
$B_\C$ and $X_\C$ the scalar extensions $B\otimes_\Q\C$ and
$X\times_{\Spec\Q}\Spec\C$, respectively.

When $X$ is of finite type over $\Q$, we denote by $X_\h$ the
analytic space associated to the algebraic $\C$-scheme $X_\C$
\cite[Chap.1,\S6]{Hartshorne}.
For a coherent sheaf $F$ on $X$, the pulled-back coherent sheaf
on $X_\C$ and the associated analytic coherent sheaf on $X_\h$
are denoted by $F_\C$ and $F_\h$, respectively.
Let $f\cl F\ra G$ be a differential operator of order one over
$\Q$.
Then it is easy to see that the $\C$-homomorphism
$f_\C\cl F_\C\ra G_\C$ on $X_\C$ obtained by scalar extension
is a differential operator of order one over $\C$.
By \cite[16.8]{EGA4}, the differential operator $f_\C\cl F_\C\ra G_\C$
is decomposed uniquly to $u\cdot d_{X_\C}^1$ where
$d_{X_\C}^1\cl F\ra{\cal P}_{X_\C}^1(F)$ is a canonical
$\C$-homomorphism \cite[16.7.5]{EGA4} and $u\cl{\cal P}_{X_\C}^1(F)\ra G$
is a $\calO_{X_\C}$-homomorphism.
For the definition of ${\cal P}_{X_\C}^1(F)$, see
\cite{EGA4}.
Since these homomorphisms are canonically pulled-back to
$X_\h$, we get a differential operator $f_\h\cl F_\h\ra G_\h$
of order one of the analytic coherent sheaves on $X_\h$.

Let $\rho$ be a cone in $N_\R$.
By the notation of the scalar extensions,
\begin{equation}
S(\rho)_\C = S(\rho)\otimes_\Q\C =\C[M\cap\rho^\vee]
\end{equation}
and
\begin{equation}
S(\rho;\sigma)_\C = S(\rho;\sigma)\otimes_\Q\C =
\C[M[\sigma]\cap\rho^\vee]
\end{equation}
for each $\sigma\in F(\rho)$.
Since $S(\rho;\sigma) = S(\rho)/P(\rho;\sigma)$, $S(\rho;\sigma)_\C$
is the quotient of $S(\rho)_\C$ by the prime ideal
$P(\rho;\sigma)_\C =\C[M\cap(\rho^\vee\setminus\sigma^\bot)]$.

We fix a finite fan $\Delta$ of $N_\R$ in this section.

We denote simply by $Z$ the toric variety $Z(\Delta)$.
Since $Z$ is normal, so are the toric variety $Z_\C = Z(\Delta)_\C$
over $\C$ and the analytic space $Z_\h = Z(\Delta)_\h$.
For each $\sigma\in\Delta$, a $T$-invariant irreducible closed
subvariety $X(\sigma)$ of $Z$ was defined in Section~3.
Hence $X(\sigma)_\C$ and $X(\sigma)_\h$ are irreducible
closed subvarieties of $Z_\C$ and $Z_\h$, respectively.

Let $n$ be an element of $N$.
The group homomorphism  $M\ra\Z$ defined by
$m\mapsto\lan m, n\ran$ induces a homomorphism of group rings
$\C[M]\ra\C[t, t^{-1}]$ and the associated morphism
$\lambda_n\cl\Spec\C[t, t^{-1}]\ra T_\C =\Spec\C[M]$, where
$t$ is the monomial corresponding to $1\in\Z$.
We call $\lambda_n$ the {\em one-parameter subgroup} associated
to $n$ \cite[1.2]{Oda1}.
We denote also by $\lambda_n$ the associated map
$\C^*\ra T_\h$ of complex Lie groups.
If $n\in N\cap\rho$ for a cone $\rho\in\Delta$, then
$\C[M\cap\rho^\vee]$ is mapped to $\C[t]$.
Hence the one-parameter subgroup is extended to a regular
map $\lambda_n\cl\C\ra Z(F(\rho))_\C\subset Z_\C$, uniquely.

For each $m\in M$, the monomial $\e(m)\in\C[M]$ is regarded as a
character $T_\h\ra\C^*$.
The composite $\e(m)\cdot\lambda_n\cl\C^*\ra\C^*$
of $\lambda_n$ and $\e(m)$ is equal to the map
$t\mapsto t^{\lan m, n\ran}$.

Since $Z_\C$ is a toric variety, the group $T_\h$ acts on
$Z_\h$ analytically.
For $a\in T_\h$ and $x\in Z_\h$, we denote by $ax$ the
corresponding point of $Z_\h$ by the action.

Let $f$ be a complex analytic function on an open subset $U$ of $Z_\h$.
For each $n\in N$, the derivation $\partial_nf$ of $f$ is defined
by
\begin{equation}
\partial_nf(x) :=\left.\frac{d}{dt}\right|_{t=1}f(\lambda_n(t)x) =
\lim_{t\ra 1}\frac{f(x) - f(\lambda_n(t)x)}{1 - t}
\end{equation}
for $x\in U$.
Since $f(\lambda_n(t)x)$ is analytic in the variables $x, t$, the
function $\partial_nf$ is anlytic on $U$.
Hence we get a $\C$-derivation
$\partial_n\cl\calO_{Z_\h}\ra\calO_{Z_\h}$ of the structure sheaf.

For each $\sigma\in\Delta$, we denote by $\calP(\sigma)_\h$
the ideal sheaf of $\calO_{Z_\h}$ defining $X(\sigma)_\h\subset Z_\h$.
If $f$ is in $\calP(\sigma)_\h(U)$, then so is $\partial_nf$
since $f$ is zero on $X(\sigma)\cap U$ and $X(\sigma)_\h$ is
closed by the action of $\lambda_n$.

%
%

\begin{Lem}
\label{Lem 4.1}
Let $\sigma$ be an element of $\Delta$ and let $y$
be a point in $X(\sigma)_\h$.
If $n$ is in $N\cap\relint\sigma$, then the endomorphism
$(\partial_n)_y$ of the stalk $(\calP(\sigma)_\h)_y$
is an automorphism as a $\C$-vector space.
\end{Lem}

\Proof
Let $U\subset Z_\h$ be an open neighborhood of $y$ such that
$\lambda_n(t)U\subset U$ for every $t$ with $|t|\leq 1$.
For $g\in\calP(\sigma)_\h(U)$, we define
\begin{equation}
f(x) :=\int_0^1\frac{g(\lambda_n(s)x)}{s}ds\;,
\end{equation}
where the integration is taken on the real interval $[0,1]$.
Since $\lambda_n(0)x$ is in $X(\sigma)_\h$, the analytic fuction
$g(\lambda_n(s)x)$ on $U\times\{s\scl|s|\leq 1\}$ has zero at the
divisor $(s = 0)$.
Hence $g(\lambda_n(s)x)/s$ is an analytic function.
Hence the integral $f$ is an analytic function on $U$.
By the definition, we have
\begin{equation}
f(\lambda_n(t)x) =\int_0^1\frac{g(\lambda_n(st)x)}{s}ds
=\int_0^t\frac{g(\lambda_n(u)x)}{u}du\;.
\end{equation}
Hence $\partial_nf = g$.
This implies that $\partial_n\cl\calP(\sigma)_\h(U)\ra\calP(\sigma)_\h(U)$
is surjective.

For $f\in\calP(\sigma)_\h(U)$, suppose that $\partial_nf = 0$.
Then $(d/dt)f(\lambda_n(t)x) = 0$ for $t\in[0, 1]$ and
we have $f(x) = f(\lambda_n(0)x) = 0$.
Since $\partial_n$ is $\C$-linear, it is also injective.

Since $U$'s with this property form a fundamental system of
neighborhood of $y$, $\partial_n\cl\calP(\sigma)_\h\ra\calP(\sigma)_\h$
is isomorphic at the stalk of $y$.
\QED

Let $D = D(\Delta)$ be the complementary reduced divisor
of $T$ in $Z$.
We set
\begin{equation}
\Omega_{Z_\h}(\log D_\h)^1 :=
\Omega_Z(\log D)_\h^1 =\calO_{Z_\h}\otimes_\Z M\;.
\end{equation}
We define the $\C$-derivation
$\partial\cl\calO_{Z_\h}\ra\Omega_{Z_\h}(\log D_\h)^1$ as follows.

Let $n, n'$ be elements of $N$.
The equality $\lambda_{n + n'}(t) =\lambda_n(t)\lambda_{n'}(t)$
holds for $t\in\C^*$.
For an analytic function $f$ on an open subset $U$ of $Z_\h$,
we have
\begin{eqnarray*}
\lefteqn{\partial_{n + n'}f(x)} \\
 & = &\lim_{t\ra 1}\frac{f(x) - f(\lambda_n(t)\lambda_{n'}(t)x)}{1 - t} \\
 & = &\lim_{t\ra 1}\frac%
{f(x) - f(\lambda_n(t)x)}{1 - t} + \lim_{t\ra 1}\frac%
{f(\lambda_n(t)x) - f(\lambda_{n'}(t)\lambda_n(t)x)}{1 - t} \\
 & = &\partial_nf(x) +\partial_{n'}f(x)
\end{eqnarray*}
for any $x\in U$.
Hence the map $n\mapsto\partial_nf\in\calO_{Z_\h}(U)$ is a
homomorphism, i.e., an element of $\Hom_\Z(N,\calO_{Z_\h}(U))$.
We define $\partial f$ to be the corresponding element of
$\Omega_{Z_\h}(\log D_h)(U)^1 =\calO_{Z_\h}(U)\otimes_\Z M$.

Let $m$ and $n$ be elements of $M$ and $N$, respectively.
Since the character $\e(m)\cl T_\h\ra\C^*$ is a homomorphism,
\begin{equation}
\e(m)(\lambda_n(t)x)
= \e(m)(\lambda_n(t))\cdot \e(m)(x)
= t^{\lan m, n\ran}\cdot \e(m)(x)\;.
\end{equation}
Hence $\partial_n\e(m) =\lan m, n\ran \e(m)$ for every $n\in N$.
This implies $\partial\e(m) =\e(m)\otimes m$.
Hence we denote the global section $1\otimes m$ of
$\Omega_{Z_\h}(\log D_h)^1$ by $\partial\e(m)/\e(m)$ for
every $m\in M$, however it is common to write it by
$d\e(m)/\e(m)$.

Let $\{m_1,\cdots, m_r\}$ be a $\Z$-basis of $M$ and let
$\{n_1,\cdots, n_r\}$ be the dual basis of $N$.
If we write $\partial f =\sum a_i\otimes m_i$, we have
$a_i =\lan\partial f, n_i\ran =\partial_{n_i}f$ for each $i$.
Hence
\begin{equation}
\label{eq df}
\partial f =\sum_{i=1}^r\partial_{n_i}f\,\frac{\partial\e(m_i)}{\e(m_i)}\;.
\end{equation}
In particular $\partial$ is a $\C$-derivation.

Since $\e(m)$'s form a $\C$-basis of the coordinate ring of
each affine toric variety, we know that this $\C$-derivation
$\partial$ is compatible with the algebraic $\C$-derivation
$\partial\cl\calO_{Z_\C}\ra\Omega_{Z_\C}(\log D_\C)^1$.

For each $0\leq i\leq r$ we set
\begin{equation}
\Omega_{Z_\h}(\log D_\h)^i :=
\bigwedge^i\Omega_{Z_\h}(\log D_\h)^1 =
\calO_{Z_\h}\otimes_\Z\bigwedge^i M\;.
\end{equation}
Since $\Omega_{Z_\h}(\log D_\h)^1$ is a free $\calO_{Z_\h}$-module
of rank $r$, $\Omega_{Z_\h}(\log D_\h)^i$ is free of rank ${}_rC_i$.
For $0 < i < r$, we define a pairing
\begin{equation}
\calO_{Z_\h}(U)\times\bigwedge^i M\lra
\Omega_{Z_\h}(\log D_\h)^{i+1}(U)
\end{equation}
by $(f, w)\mapsto\partial f\wedge w$ which induces a $\C$-homomorphism
\begin{equation}
\partial\cl\Omega_{Z_\h}(\log D_\h)^i\ra\Omega_{Z_\h}(\log D_\h)^{i+1}
\end{equation}
of sheaves which we denote also by $\partial$.
We check easily that $\partial\cdot \partial = 0$, and we get a
$\partial$-complex
$\Omega_{Z_\h}(\log D_\h)^\bullet$ which we call the
{\em logarithmic de Rham complex} on $Z_\h$.
For any $\sigma\in\Delta$, we see easily by the description
(\ref{eq df}) of $\partial f$ that
$\calP(\sigma)_\h\Omega_{Z_\h}(\log D_\h)^\bullet$
is a subcomplex of $\Omega_{Z_\h}(\log D_\h)^\bullet$.

For each $\sigma\in\Delta$, we denote by $\bar I_\sigma$  the closed
immersion $X(\sigma)_\h\ra Z_\h$ and by $I_\sigma$ the immersion
$T[\sigma]_\h\ra Z_\h$.
For a $\C$-sheaf $F$ on $Z_\h$, we denote by $\bar I_\sigma^*F$
and $I_\sigma^*F$ the pull-back of $F$ to $X(\sigma)_\h$ and
$T[\sigma]_\h$, respectively.
Note that, even if $F$ is an $\calO_{Z_\h}$-module, the pull-back
is taken as a $\C$-sheaf.

%
%

\begin{Lem}
\label{Lem 4.2}
For every $\sigma\in\Delta$, the $\partial$-complex
$\bar I_\sigma^*(\calP(\sigma)_\h\Omega_{Z_\h}(\log D_\h))^\bullet$
on $X(\sigma)_\h$ is homotopically equivalent to the zero complex.
\end{Lem}

\Proof
We have to show that the identity map of this $\partial$-complex is
homotopic to the zero map.
Take an element $n_0$ in $N\cap\relint\sigma$.
Let $h$ be the $\calO_{Z_\h}$-homomorphism of degree $-1$ of the
graded $\calO_{Z_\h}$-module $\Omega_{Z_\h}(\log D_\h)^\bullet$
induced by the right interior product
$i(n_0)\cl\bigwedge^\bullet M\ra\bigwedge^\bullet M$.
Then the $\C$-homomorphism $h\cdot \partial + \partial\cdot h$ on
$\Omega_{Z_\h}(\log D_\h)^\bullet = \calO_{Z_\h}\otimes_\Z\bigwedge^\bullet M$
is equal to $\partial_{n_0}\otimes 1$ (cf. \cite[13.4]{Danilov}).
By Lemma~\ref{Lem 4.1}, this induces an automorphism of
$\bar I_\sigma^*(\calP(\sigma)_\h\Omega_{Z_\h}(\log D_\h))^\bullet$
as a $\C$-sheaf.
Let $u$ be the inverse isomorphism of
$\bar I_\sigma^*(h\cdot \partial + \partial\cdot h)$.
Since
\begin{equation}
(h\cdot \partial + \partial\cdot h)\cdot\partial =
\partial\cdot(h\cdot \partial + \partial\cdot h) =
\partial\cdot h\cdot \partial\;,
\end{equation}
we have
$\bar I_\sigma^*\partial\cdot u = u\cdot\bar I_\sigma^*\partial
=u\cdot I_\sigma^*(\partial\cdot h\cdot \partial)\cdot u$.
Set $h' := u\cdot\bar I_\sigma^*h$.
Then we have
\begin{equation}
h'\cdot\bar I_\sigma^*\partial +\bar I_\sigma^*\partial\cdot h' =
u\cdot I_\sigma^*(h\cdot \partial + \partial\cdot h) = 1\;.
\end{equation}
\QED

Let $\sigma$ be an element of $\Delta$.
We set $N[\sigma] := N/(N\cap(\sigma+(-\sigma)))$.
This is naturally the dual $\Z$-module of
$M[\sigma] = M\cap\sigma^\bot$.
It is known that the closed subscheme $X(\sigma)$ of $Z$ is
naturally identified with the toric vatiety $Z(\Delta[\sigma])$
with the torus $T[\sigma]$ \cite[Cor.1.7]{Oda1}.
Furthermore, the action of $T[\sigma]$ on $Z(\Delta[\sigma])$
is equivariant with that of $T$ with respect to the natural
surjection $T\ra T[\sigma]$.

We set $D(\sigma)_\h := X(\sigma)_\h\setminus T[\sigma]_\h$
and
\begin{equation}
\Omega_{X(\sigma)_\h}(\log D(\sigma)_\h)^\bullet :=
\calO_{X(\sigma)_\h}\otimes_\Z{\bigwedge}\!{}^\bullet M[\sigma]\;.
\end{equation}
Then $\Omega_{X(\sigma)_\h}(\log D(\sigma)_\h)^\bullet$ has the
$\partial$-complex structure so that it is identified with
the logarithmic de Rham complex of $Z(\Delta[\sigma])_\h$
similarly as $\Omega_{Z_\h}(\log D_\h)^\bullet$ for
$Z_\h =Z(\Delta)_\h$.

%
%

\begin{Lem}
\label{Lem 4.3}
Let $\sigma$ be an elemet of $\Delta$.
For $\rho\in\Delta$ with $\sigma\prec\rho$, the $\partial$-complex
$\bar I_\rho^*(\calP(\rho)_\h\Omega_{X(\sigma)_\h}(\log D(\sigma)_\h))^\bullet$
of $\C$-sheaves on $X(\rho)_\h$ is homotopically equivalent
to the zero complex.
\end{Lem}

\Proof
For the toric variety $X(\sigma)_\h = Z(\Delta[\sigma])_\h$, $X(\rho)_\h$
is the closed subvariety associated to $\rho[\sigma]\in\Delta[\sigma]$.
The closed subvariety $X(\rho)_\h$ of $X(\sigma)_\h$ is defined by the
image of $\calP(\rho)_\h$ in $\calO_{X(\sigma)_\h}$.
Hence, this is a consequence of Lemma~\ref{Lem 4.2} applied
to the toric variety $Z(\Delta[\sigma])_\h$.
\QED

Let $\rho$ be an element of $\Delta$ and $V$ an object of $\GM(A(\rho))$.
Then the $\partial$-complex $\Lambda_Z^\rho(V)^\bullet$ defined in
Section~3 induces a $\partial$-complex $\Lambda_{Z_\C}^\rho(V)^\bullet$
and its analytic version $\Lambda_{Z_\h}^\rho(V)^\bullet$.
By Proposition~\ref{Prop 3.1}, $\Lambda_{Z_\h}^\rho(V)^\bullet$ is
isomorphic to a finite direct sum of dimension shifts of
$\Omega_{X(\rho)_\h}(\log D(\rho)_\h))^\bullet$.
In particular, we get the following corollary by Lemma~\ref{Lem 4.3}.

%
%

\begin{Cor}
\label{Cor 4.4}
Let $\sigma$ be an element of $\Delta$ and $V$ in $\GM(A(\sigma))$.
For $\rho\in\Delta$ with $\sigma\prec\rho$, the $\partial$-complex
$\bar I_\rho^*(\calP(\rho)_\h\Lambda_{Z_h}^\sigma(V))^\bullet$
on $X(\rho)_\h$ is homotopically equivalent to the zero complex.
\end{Cor}

We recall some general notation of derived categories
(cf.\cite{Verdier}).

Let $X$ be a locally compact topological space and let
$A(\C_X)$ be the abelian category of the $\C$-sheaves on $X$.
We denote by $C^+(\C_X)$ and $D^+(\C_X)$ the category of the
complexes bounded below in $A(\C_X)$ and the derived category
of it, respectively.
For a continuous map $f\cl X\ra Y$, let
$f^*\cl A(\C_Y)\ra A(\C_X)$ be the the pull-back functor
which is exact.
We denote also by $f^*$ the induced functors
$C^+(\C_Y)\ra C^+(\C_X)$ and $D^+(\C_Y)\ra D^+(\C_X)$.

If the direct image functor with proper support
$f_!\cl A(\C_X)\ra A(\C_Y)$ is of finite cohomological dimension,
the functor $f^!\cl D^+(\C_Y)\ra D^+(\C_X)$ is defined
\cite[2.2]{Verdier}.
Note that this condition is satisfied for any regular morphisms
of finite dimensional analytic spaces.

For $V^\bullet\in\CGM(A(\sigma))$, the bicomplex
$\Lambda_{Z_\h}^\sigma(V)^{\bullet,\bullet}$ is defined
as the analytic version of $\Lambda_Z^\sigma(V)^{\bullet,\bullet}$.
The associated single complex is denoted by
$\Lambda_{Z_\h}^\sigma(V)^\bullet$.

For $L^\bullet\in\CGEM(\Delta)$, we get the bicomplex
$\Lambda_{Z_\h}(L)^{\bullet,\bullet}$ and its associated
single complex $\Lambda_{Z_\h}(L)^\bullet$, similarly.
We denote by $\delta$ the coboundary map of
$\Lambda_{Z_\h}(L)^\bullet$.

%
%

\begin{Prop}
\label{Prop 4.5}
For $\rho\in\Delta$ and a $d$-complex $L^\bullet$ in
$\CGEM(\Delta)$, there exist a quasi-isomorphism
\begin{equation}
I_\rho^*\Lambda_{Z_\h}(L)^\bullet\simeq
I_\rho^*\Lambda_{Z_\h}^\rho(\i_\rho^*L)^\bullet
\end{equation}
as $\delta$-complexes of $\C$-sheaves on $T[\rho]_\h$ and an
isomorphism
\begin{equation}
I_\rho^!\Lambda_{Z_\h}(L)^\bullet\simeq
I_\rho^*\Lambda_{Z_\h}^\rho(\i_\rho^!L)^\bullet
\end{equation}
in the derived category $D^+(\C_{T[\rho]_\h})$.
\end{Prop}

\Proof
Let $\sigma$ be an element of $F(\pi)$ with $\sigma\prec\rho$.
For $V$ in $\GM(A(\sigma))$, there exists an exact sequence
\begin{equation}
0\lra\calP(\rho)_\h\Lambda_{Z_\h}^\sigma(V)^\bullet\lra
\Lambda_{Z_\h}^\sigma(V)^\bullet\mathop{\lra}\limits^{\lambda(V)}
\Lambda_{Z_\h}^\rho(V_{A(\rho)})^\bullet\lra 0\;.
\end{equation}
In particular, the homomorphism $I_\rho^*\lambda(V)$ of
$\partial$-complexes is a quasi-isomorphism by Corollary~\ref{Cor 4.4}.

For each $L^i$, we get a quasi-isomorphism
\begin{equation}
I_\rho^*\lambda(L^i)\cl I_\rho^*\Lambda_{Z_\h}(L^i)^\bullet\lra
I_\rho^*\Lambda_{Z_\h}^\rho(\i_\rho^*L^i)^\bullet
\end{equation}
as a collection of quasi-isomorphisms $I_\rho^*\lambda(L^i(\sigma))$
for $\sigma\in F(\rho)$.

We get the first quasi-isomorphism as the collection of
$I_\rho^*\lambda(L^i)$'s for $i\in\Z$.

Since $\i_\rho^!L^\bullet = L(\rho)^\bullet$, it is enough to
show that $\bar I_\rho^!\Lambda_{Z_\h}^\sigma(L^i(\sigma))$'s are
quasi-isomorphic to the zero complex for all
$\sigma\in F(\rho)\setminus\{\rho\}$ and $i\in\Z$.
By proposition~\ref{Prop 3.1}, $\Lambda_{Z_\h}^\sigma(L^i(\sigma))$
is isomorphic to the finite direct sum of dimension shifts
of the logarithmic de Rham complex
$\Omega_{X(\sigma)_\h}(\log D(\sigma)_\h)^\bullet$ on $X(\sigma)_\h$.
By \cite[Prop.1.2]{Oda2}, the logarithmic de Rham complex
is quasi-isomorphic to the direct image $Rj_*\C_{T[\sigma]}$ for
the open immersion $j\cl T[\sigma]_\h\ra X(\sigma)_\h$.
Since $T[\sigma]_\h\cap X(\rho)_\h =\emptyset$,
$\bar I_\rho^!Rj_*\C_{T[\sigma]_\h}$ is equivalent to zero.
Hence $\bar I_\rho^!\Lambda_{Z_\h}^\sigma(L^i(\sigma))$ is also zero
in the derived category.
\QED

For any $\Q$-vector space $W$, we denote by $\bar W$ the scalar
extension $W\otimes_\Q\C$.
If $W$ is graded, so is $\bar W$.

Let $\rho$ be an element of $\Delta$.
For $V$ in $\GM(A(\rho))$, we denote by $\bar V_{T[\rho]_\h}$ the
constant sheaf on $T[\rho]_\h$ with the stalk $\bar V$.
We regard it as a $\partial$-complex of $\C$-sheaves by setting
$\bar V_{T[\rho]_\h}^i := (\bar V_i)_{T[\rho]_\h}$ and $\partial = 0$.

For $N[\rho] := N/N(\rho)$, we set $\det[\rho]_\Q :=\det N[\rho]_\Q$.
We denote by $\Det[\rho]$ the graded $\Z$-module defined by
$(\Det[\rho])_{-r+r_\rho} :=\det[\rho]$ and
$(\Det[\rho])_i :=\{0\}$ for $i\not= -r+r_\rho$.
Here note that $\rank N[\rho] = r-r_\rho$.

We define a homomorphism
$(\bar V\otimes\Det[\rho])_{T[\rho]_\h}^\bullet\ra
I_\rho^*(\Lambda_{Z_\h}^\rho(V))^\bullet$ as follows.

We take a $\Q$-linear subspace $H$ of $N_\Q$ such that
$N_\Q = N(\rho)_\Q\oplus H$.
Then we have a natural isomorphisms
$N[\rho]_\Q = N_\Q/N(\rho)_\Q\simeq H$ and
$\det[\rho]_\Q =\det N[\rho]_\Q\simeq\det H\subset A_{-r+r_\rho}$.
We denote by $\Det H$ the graded $\Q$-vector space defined by
$(\Det H)_{-r+r_\rho} :=\det H$ and $(\Det H)_i :=\{0\}$ for
$i\not= -r+r_\rho$.
By Lemma~\ref{Lem 1.3}, $V_A = V\otimes_{A(\rho)}A$ is equal
to the free $A^*[\rho]$-module
$(V\otimes_\Q\Det H)\otimes_\Q A^*[\rho]$.

The restriction of $\partial$ of
$\Lambda_Z^\rho(V)^\bullet =\calO_Z\otimes_\Q V_A$ to the
constant subsheaf $\Q\e(0)\otimes_\Q V_A$ is zero.
Hence the composite of the natural homomorphisms
\begin{equation}
\bar V\otimes\Det[\rho]\lra\bar V\otimes_\Q\Det H\lra
\bar V_A\lra \C\e(0)\otimes_\C\bar V_A
\end{equation}
defines a homomorphism
$(\bar V\otimes\Det[\rho])_{T[\rho]_\C}\ra%
\Lambda_{Z_\C}^\rho(V)|_{T[\rho]_\C}$
of $\partial$-complexes of $\C$-sheaves on the algebraic torus
$T[\rho]_\C$

We denote by $\phi_V^H$ the associated homomorphism
\begin{equation}
(\bar V\otimes\Det[\rho])_{T[\rho]_\h}^\bullet\lra
I_\rho^*(\Lambda_{Z_\h}^\rho(V))^\bullet
\end{equation}
on the smooth analytic space $T[\rho]_\h$.

%
%

\begin{Prop}
\label{Prop 4.6}
Let $\rho\in\Delta$ and let $V$ be an object of $\GM(A(\rho))$.
Then the homomorphism $\phi_V^H$ is a quasi-isomorphism for any
$H\subset N_\Q$ with $N_\Q = N(\rho)_\Q\oplus H$.
For another $K\subset N_\Q$ with $N_\Q = N(\rho)_\Q\oplus K$,
the homomorphisms $\phi_V^H$ and $\phi_V^K$ are locally homotopic.
\end{Prop}

\Proof
Since $V$ is of finite dimension, we prove the proposition
by induction on the dimension of $V$.
If $V =\{0\}$, then the $\partial$-complexes are trivial and the
assertion is clear.
Assume $\dim V\geq 1$.
Let $k$ be the maxiaml integer with $V_k\not=\{0\}$.
We take a vector subspace $V'_k\subset V_k$ of codimension one.
By setting $V'_i := V_i$ for $i\not= k$, we get a homogeneous
subspace $V'\subset V$ of codimension one.
Since $A(\rho)$ is graded negatively, $V'$ is an object of $\GM(A(\rho))$.
We get an exact sequence
\begin{equation}
0\lra V'\lra V\lra\Q(-k)\lra 0
\end{equation}
of graded $A(\rho)$-modules, where $Q(a)$ is the graded
$\Q$-vector space defined by $Q(a)_{-a} :=\Q$ and
$Q(a)_i :=\{0\}$ for $i\not= -a$.
Set $c := r - r_\rho$.
Then we get a commutative diagram
\begin{equation}
\begin{array}{ccccccccc}
0 &\ra   &(\bar V'\otimes\Det[\rho])_{T[\rho]_\h}^\bullet &\ra
&(\bar V\otimes\Det[\rho])_{T[\rho]_\h}^\bullet &\ra
&(\C({-}k)\otimes\Det[\rho])_{T[\rho]_\h}^\bullet &\ra & 0 \\
 &   &\hbox{  }\downarrow\phi_{V'}^H &   &\hbox{  }\downarrow\phi_V^H &   &
\hbox{  }\downarrow\phi_{\Q(-k)}^H &   & \\
0 &\ra   & I_\rho^*(\Lambda_{Z_\h}^\rho(V'))^\bullet &\ra
&I_\rho^*(\Lambda_{Z_\h}^\rho(V))^\bullet &\ra
&I_\rho^*(\Lambda_{Z_\h}^\rho(\Q(-k)))^\bullet &\ra & 0
\end{array}
\end{equation}
By the induction assumption, $\phi_{V'}^H$ is a quasi-isomorphism
of $\partial$-complexes of $\C$-sheaves.
On the other hand, $I_\rho^*(\Lambda_{Z_\h}^\rho(\Q(-c)))^\bullet$ is equal
to the analytic de Rham complex on the complex manifold $T[\rho]_\h$.
Since $\phi_{\Q(-k)}^H$ is the dimension shift of the natural
homomorphism $\C_{T[\rho]_\h}\lra\Omega_{T[\rho]_\h}^\bullet$, this is a
quasi-isomorphism, by the complex analytic version of the
Poincar\'e Lemma.
Hence $\phi_V^H$ is also a quasi-isomorphism.

Clearly, $\phi_{\Q(-k)}^H$ does not depend on the choice of $H$.
Hence $\phi_{\Q(-k)}^H -\phi_{\Q(-k)}^K = 0$.
By assumption, $\phi_{V'}^H -\phi_{V'}^K$ is locally homotopic
to zero.
Hence, we know that $\phi_V^H -\phi_V^K$ gives zero maps on the
cohomology sheaves.
Since $\bar V_{T[\rho]_\h}$ is a locally free $\C$-sheaf,
$\phi_V^H -\phi_V^K$ is locally homotopic to zero.
\QED

Recall that $Z_\h$ has the docomposition
$\bigcup_{\sigma\in\Delta}T[\sigma]_\h$ into $T_\h$-orbits.
For each integer $0\leq i\leq r$, we define
\begin{equation}
Z_\h^{2i} = Z_\h^{2i+1} =:
\bigcup_{\scriptstyle\sigma\in\Delta\atop\scriptstyle r_\sigma\geq r-i}
T[\sigma]_\h\;.
\end{equation}
Then we have a filtration
\begin{equation}
Z_\h = Z_\h^{2r}\supset Z_\h^{2r-1}\supset\cdots
\supset Z_\h^2\supset Z_\h^1\supset Z_\h^0
\end{equation}
of $Z_\h$ satisfying the conditions of the topological
stratification in \cite[1.1]{GM2}.

The intersection complex of a stratified space is defined for a
sequence of integers $(\p(2),\p(3),\p(4),\cdots)$ which
is called a {\em perversity} \cite[2.0]{GM2}.
Since $Z_\h$ is a complex analytic space of dimension $r$, only
$\p(i)$'s for even $i$ less than or equal to $2r$ are relevant
for the intersection complex.

Let $\p = (\p(2),\p(4),\p(6),\cdots,\p(2r))$ be a sequence of
integers with $\p(2) = 0$ and $\p(2i)\leq\p(2i+2)\leq\p(2i)+2$ for
$i=1,\cdots, r-1$ as a perversity for $Z_\h$.
We denote by the same symbol $\p$ the perversity on $\Delta$
defined by $\p(\sigma) :=\p(2r_\sigma) -r_\sigma + 1$
for $\sigma\in\Delta\setminus\{\Zero\}$.
Note that, for the middle perversity $\m := (0, 1, 2,\cdots, r-1)$,
we have $\m(\sigma) = 0$ for all $\sigma\in\Delta\setminus\{\Zero\}$.

We consider the $\delta$-complex $\Lambda_{Z_\C}(\ic_\p(\Delta))^\bullet$
which is the associated single complex of the bicomplex
$\Lambda_Z(\ic_\p(\Delta))^{\bullet,\bullet}$.

By \cite[3.3, AX1]{GM2}, the intersection complex $\IC_\p(Z_\h)^\bullet$
in the derived category $D^{\rm b}(\C_{Z_\h})$ of bounded
complexes of $\C$-sheaves is characterized by the following
properties.

For the convenience of our use, we set
$F :=\IC_\p(Z_\h)[-r]^\bullet$.
Note that $n$ in \cite[3.3]{GM2} is $2r$ in our case.

(a) The restriction of $F$ to $T_\h$ is quasi-isomorphic
to $\C_{T_\h}[r]$.

(b) $\calH^i(F)$ is a trivial sheaf for $i < -r$.

(c) For any $\sigma\in\Delta\setminus\{\Zero\}$,
$\H^i(F_x) =\{0\}$ for $i \geq\p(2r_\sigma) - r + 1$.

(d) For any $\sigma\in\Delta\setminus\{\Zero\}$,
$\H^i((I_\sigma^!F)_x) =\{0\}$ for $i \leq\p(2r_\sigma) - r + 1$.

%
%

\begin{Thm}
\label{Thm 4.7}
Let $\p = (\p(2),\p(4),\p(6),\cdots,\p(2r))$ be a perversity.
The complex of $\C$-sheaves $\Lambda_{Z_\h}(\ic_\p(\Delta))^\bullet$
is isomorphic to $\IC_\p(Z_\h)[-r]^\bullet$.
\end{Thm}

\Proof
We set $L^\bullet :=\Lambda_{Z_\h}(\ic_\p(\Delta))^\bullet$.
It is sufficient to show that $L^\bullet$ satisfies the above
properties.
Let $\sigma$ be in $\Delta\setminus\{\Zero\}$.
Then the stratum $T[\sigma]_\h$ is of dimension $2r - 2r_\sigma$
in the real dimension.
By Propositions~\ref{Prop 4.5} and \ref{Prop 4.6}, we have an
isomorphism
\begin{equation}
\H^i(L^\bullet_x)\simeq
\H^{i+r-r_\sigma}(\i_\sigma^*\ic_\p(\Delta)^\bullet)\otimes_\Q\C
\end{equation}
for every point $x\in T[\sigma]_\h$.
By the condition (3) of Theorem~\ref{Thm 2.9}, the cohomology
$\H^i(\i_\sigma^*\ic_\p(\Delta)^\bullet)$ vanishes for
$i\geq\p(\sigma) =\p(2r_\sigma) - r_\sigma + 1$.
Hence $\H^i(L_x^\bullet) = 0$ for
$i\leq\p(2r_\sigma) - r + 1$.
By Propositions~\ref{Prop 4.5} and \ref{Prop 4.6}, we also have
\begin{equation}
\H^i((I_\sigma^!L)_x^\bullet)\simeq
\H^{i+r-r_\sigma}(\i_\sigma^!\ic_\p(\Delta)^\bullet)\otimes_\Q\C
\end{equation}
for all $x\in T[\sigma]_\h$.
The condition (2) of Theorem~\ref{Thm 2.9} implies
$\H^i(\i_\sigma^!\ic_\p(\Delta)^\bullet) = 0$ for
$i\leq\p(\sigma) =\p(2r_\sigma) - (r_\sigma) + 1$.
Hence $\H^i((I_\sigma^!L)_x^\bullet) = 0$ for $i\leq\p(2r_\sigma) - r + 1$.
Hence, by the Axioms [AX1] of \cite[3.3]{GM2}, $L^\bullet$ is the
intersection complex with the perversity $\p$.
\QED

%
%

\begin{Prop}
\label{Prop 4.8}
Assume that $\Delta$ is a complete fan.
Let $L^\bullet$ be an object of $\CGEM(\Delta)$.
Then there exists a natural isomorphism
\begin{equation}
\H^k(\Lambda_{Z_\h}(L)^\bullet)\simeq
\bigoplus_{p+q=k}\H^p(\Gamma(L)^\bullet)_q\otimes_\Q\C
\end{equation}
for every integer $k$.
\end{Prop}

\Proof
Since any $X(\sigma)$ is a complete toric variety,
$\H^i(X(\sigma),\calO_{X(\sigma)}) = 0$
for every $i > 0$ \cite[Cor.2.8]{Oda1}.
By \cite{GAGA}, we have
\begin{equation}
\H^i(X(\sigma)_\h,\calO_{X(\sigma)_\h}) =\left\{
\begin{array}{cc}
  \C    & \;\;\;\hbox{ if }\;\;  i = 0 \\
  \{0\} & \;\;\;\hbox{ if }\;\;  i > 0
\end{array}
\right.
\end{equation}

Since $\Lambda_{Z_\h}^\sigma(L^i(\sigma))$ is a free
$O_{X(\sigma)_\h}$-module for any $\sigma$ and $i$,
we have an isomorphism
\begin{equation}
\H^k(\Lambda_{Z_\h}(L)^\bullet) =
\H^k(\Gamma(Z_\h,\Lambda_{Z_\h}(L))^\bullet)\;,
\end{equation}
where $\Gamma(Z_\h,\Lambda_{Z_\h}(L))^\bullet$ is the single
complex of $\C$-vector spaces associated to the bicomplex
$\Gamma(Z_\h,\Lambda_{Z_\h}(L))^{\bullet,\bullet}$.
Since the global sections are locally $M$-homogeneous of
degree $0\in M$, $d_2 =\partial$ of the last bicomplex is zero.
Hence the spectral sequence \begin{equation}
E_1^{p,q} :=\H^p(\Gamma(Z_\h,\Lambda_{Z_\h}(L))_q^\bullet)
\Longrightarrow
\H^{p+q}(\Gamma(Z_\h,\Lambda_{Z_\h}(L))^\bullet)
\end{equation}
degenerates at $E_1$-terms.
Consequently, the cohomology group
$\H^k(\Gamma(Z_\h,\Lambda_{Z_\h}(L))^\bullet)$ is equal to the
direct sum
\begin{equation}
\bigoplus_{q\in\Z}\H^{k-q}(\Gamma(Z_\h,\Lambda_{Z_\h}(L))_q^\bullet)\\
= \bigoplus_{q\in\Z}\H^{k-q}(\Gamma(L)^\bullet)_q\otimes_\Q\C\;.
\end{equation}
\QED


\section{The Serre duality}



\setcounter{equation}{0}

In this section, we again fix a finite fan $\Delta$ of $N_\R$
and let $Z := Z(\Delta)$ be the associated toric variety defined
over $\Q$.

Let $\p$ be a perversity on $\Delta$.
Then $\Lambda_Z(\ic_\p(\Delta))^{\bullet,\bullet}$ is a
bicomplex of coherent $\calO_Z$-modules such that $d_1 = d$
is a $\calO_Z$-homomorphism and $d_2 =\partial$ is a
differential operator of order one.
Note that $\Lambda_Z(\ic_\p(\Delta))^{i,j} =\{0\}$
unless $0\leq i\leq r$ and $-r\leq j\leq 0$ by
Proposition~\ref{Prop 2.11}.

For each integer $j$, we set
\begin{equation}
\Omega_j(\p; Z) :=\Lambda_Z(\ic_\p(\Delta))_{-j}^\bullet\;.
\end{equation}
Clearly, $\Omega_j(\p; Z)$ is the zero complex unless $0\leq j\leq r$.
Since $d$ is an $\calO_Z$-homomorphism, $\Omega_j(\p; Z)$ is a
finite complex in the category of coherent $\calO_Z$-modules.

For the top perversity $\t$, we have
\begin{equation}
\Omega_j(\t; Z)^i =
\bigoplus_{\sigma\in\Delta(i)}\det(\sigma)\otimes
\bigwedge^{-j}N[\sigma]\otimes\calO_{X(\sigma)}\;,
\end{equation}
where $N[\sigma] := N/N(\sigma)$.
For $\sigma\in\Delta(i)$ and $\tau\in\Delta(i+1)$, the
$(\sigma,\tau)$-component of the coboundary map is the
tensor product of $q'_{\sigma/\tau}\cl\det(\sigma)\ra\det(\tau)$
and the natural surjection
\begin{equation}
\bigwedge^{-j}N[\sigma]\otimes\calO_{X(\sigma)}\ra
\bigwedge^{-j}N[\tau]\otimes\calO_{X(\tau)}
\end{equation}
if $\sigma\prec\tau$ and is the zero map otherwise.

For finite complexes $B^\bullet, C^\bullet$ of coherent
sheaves on $Z$, we denote by $\calH om_{\calO_Z}(C, B)^{\bullet,\bullet}$
the bicomplex with the component
\begin{equation}
\calH om_{\calO_Z}(C, B)^{i,j} :=
\calH om_{\calO_Z}(C_j, B_i)
\end{equation}
for each pair $(i, j)$ of integers.
We denote by $\calH om_{\calO_Z}(C, B)^\bullet$ the associated
single complex.
The coboundary map
\begin{equation}
d^k\cl\calH om_{\calO_Z}(C, B)^k\lra\calH om_{\calO_Z}(C, B)^{k+1}
\end{equation}
is determined as follows.
For an integer $i$, the restriction of $d^k$ to the component
$\calH om_{\calO_Z}(C^i, B^{i+k})$ is the sum of
\begin{equation}
(d_B^{i+k})_*\cl
\calH om_{\calO_Z}(C^i, B^{i+k})\lra\calH om_{\calO_Z}(C^i, B^{i+k+1})
\end{equation}
and
\begin{equation}
(-1)^{k+1}(d_C^{i-1})^*
\calH om_{\calO_Z}(C^i, B^{i+k})\lra\calH om_{\calO_Z}(C^{i-1}, B^{i+k})
\end{equation}
(cf.\cite[1.1.10,(ii)]{Deligne})

%
%

\begin{Thm}
\label{Thm 5.1}
Let $L^\bullet$ be an object of $\CGEM(\Delta)$.
For each integer $q$, there exists a natural isomorphism
\begin{equation}
\calH om_{\calO_Z}(\Lambda_Z(L)_q,\Omega_0(\t; Z)\otimes\det N)^\bullet
\simeq\Lambda_Z(\D(L))_{-r-q}^\bullet\;.
\end{equation}
\end{Thm}

\Proof
For each pair $(i, j)$ of integers, we have
\begin{eqnarray}
\lefteqn{
\calH om_{\calO_Z}(\Lambda_Z(L)_q,
\Omega_0(\t; Z)\otimes\det N)^{i,j}} \\
 & = &
\bigoplus_{(\sigma,\rho)}\Hom_\Q((L(\sigma)_A^{-j})_q,
\det(\rho)\otimes\det N)\otimes_\Q\calO_{X(\rho)}\;,
\end{eqnarray}
where the sum is taken over all pairs $(\sigma,\rho)$ with
$\sigma\in\Delta$, $\rho\in\Delta(i)$ and $\sigma\prec\rho$.
We identify $\Hom_\Q((L(\sigma)_A^{-j})_q,\det N_\Q)$ with
$\d_N(L(\sigma)_A^{-j})_{-r-q}$ through the set of right
operations $\d_N^{\rm right}(L(\sigma)_A^{-j})_{-r-q}$.
Hence
\begin{eqnarray}
\lefteqn{
\calH om_{\calO_Z}(\Lambda_Z(L)_q,
\Omega_0(\t; Z)\otimes\det N)^k} \\
\label{Groth dual}
 & = &
\bigoplus_{\rho\in\Delta}
\det(\rho)\otimes
\d_N((\i_\rho^*L^{-r_\rho+k})_A)_{-r-q}\otimes_\Q\calO_{X(\rho)} \\
 & = &
\bigoplus_{\rho\in\Delta}
\det(\rho)\otimes
(\d_\rho(\i_\rho^*L^{-r_\rho+k})_A)_{-r-q}\otimes_\Q\calO_{X(\rho)} \\
 & = &
\Lambda_Z(\D(L))_{-r-q}^k\;.
\end{eqnarray}
The equality of the coboundary maps is also checked.
Namely, for $\rho,\mu\in\Delta$, the $(\rho,\mu)$-component
of $d^k$'s are nonzero only for (a)
$\rho\prec\mu$ and $r_\mu = r_\rho + 1$, or (b) $\rho =\mu$.
In case (a), they are both equal to the tensor product of
$q'_{\rho/\mu}$, inclusion map
\begin{equation}
\d_N((\i_\rho^*L^{-r_\rho+k})_A)_{-r-q}\ra
\d_N((\i_\mu^*L^{-r_\rho+k})_A)_{-r-q}
\end{equation}
and the natural surjection $\calO_{X(\rho)}\ra\calO_{X(\mu)}$
in the description (\ref{Groth dual}).
In case (b), they are both equal to
$(-1)^{k+1}{\rm id}\otimes\d_N(\i_\rho^*(d_L)^{k-1})\otimes{\rm id}$.
\QED

Let $S$ be a scheme of finite type over a field and let
$D_{\rm coh}^+(S)$ be the derived category of complexes
bounded below of $\calO_S$-modules with coherent cohomologies.
The Grothendieck theory of residues and duality \cite{RD} says that
there exists an object of $D_{\rm coh}^+(S)$ which is called the
dualizing complex, and the Serre duality theorem for nonsingular
projective varieties is generalized to $S$ by using the the
dualizing complex in place of the canonical invertible sheaf.

For the toric variety $Z$, the dualizing complex which is
denoted by $\calC^\bullet(Z,\Omega_0^\vee)$ in
\cite[\S5]{Ishida2} is described as follows.
For the compatibility with the notation of this paper,
we write it by $\calC(Z,\Omega_0^\vee)^\bullet$.

For each integer $-r\leq i\leq 0$, we set
\begin{equation}
\calC(Z,\Omega_0^\vee)^i :=
\bigoplus_{\sigma\in\Delta(r+i)}\calO_{X(\sigma)}\otimes\det M[\sigma] =
\bigoplus_{\sigma\in\Delta(r+i)}\Omega_{X(\sigma)}(\log D(\sigma))^{-i}\;.
\end{equation}
For $\sigma,\tau\in\Delta$, the subvariety $X(\tau)$ of $Z$
is contained in $X(\sigma)$ if and only if $\sigma\prec\tau$.
When $\sigma\prec\tau$, let
$\varphi_{\sigma/\tau}\cl\calO_{X(\sigma)}\ra\calO_{X(\tau)}$ be the
natural surjection.
The component of the coboudary map
\begin{equation}
\begin{array}{ccc}
   \calC^i(Z,\Omega_0^\vee)
 & \mathop{\lra}\limits^{\textstyle d}
 & \calC^{i+1}(Z,\Omega_0^\vee) \\
   \| & & \| \\
\bigoplus_{\sigma\in\Delta(r+i)}\calO_{X(\sigma)}\otimes\det M[\sigma]
 & &
\bigoplus_{\tau\in\Delta(r+i+1)}\calO_{X(\tau)}\otimes\det M[\tau]
\end{array}
\end{equation}
for $\sigma\in\Delta(r{+}i)$ and $\tau\in\Delta(r{+}i{+}1)$
is defined to be $\varphi_{\sigma/\tau}\otimes q_{\sigma/\tau}$
if $\sigma\prec\tau$ and the zero map otherwise, where
$\varphi_{\sigma/\tau}$ is the natural surjection
$\calO_{X(\sigma)}\ra\calO_{X(\tau)}$.
For the definition of $q_{\sigma/\tau}$, see \cite[\S1]{Ishida2}.
With respect to the identifications
$\det M[\sigma]\otimes\det N =\det(\sigma)$ and
$\det M[\tau]\otimes\det N =\det(\tau)$, the isomorphism
$q'_{\sigma/\tau}$ defined in \S2 is equal to
$q_{\sigma/\tau}\otimes 1_{\det N}$.

As a special case of \cite[Thm.5.4]{Ishida2}, the $d$-complex
$\calC(Z,\Omega_0^\vee)^\bullet$ is quasi-isomorphic to
the residual complex $f_Z^\Delta\Q$ \cite[VI,\S3]{RD} for the
structure morphism $f_Z\cl Z\ra\Spec\Q$, i.e., it is a global
dualizing complex of the $\Q$-scheme $Z$.

For a finite complex $B^\bullet$ of $\calO_Z$-modules with
coherent cohomology sheaves, the object
$R\calH om(B, f_Z^\Delta\Q)^\bullet$ in the derived category
$D_{\rm coh}^+(Z)$ is called the Grothendieck dual of
$B^\bullet$.

It is easy to see that $\Lambda_Z(\ic_\t(\Delta))_0^\bullet$ is
isomorphic to $(\calC(Z,\Omega_0^\vee))\otimes(\det N))[-r]^\bullet$.

%
%

\begin{Cor}
\label{Cor 5.2}
For each integer $0\leq j\leq r$, the Grothendieck dual of the
$d$-complex $\Lambda_Z(L)_{-j}^\bullet$ is quasi-isomorphic to
$\Lambda_Z(\D(L))_{j-r}[r]^\bullet$.
\end{Cor}

\Proof
Since $\Lambda_Z(L)_{-j}^\bullet$ is a direct sum of free
$\calO_{X(\sigma)}$-modules for $\sigma\in\Delta$, the dual
$R\calH om(\Lambda_Z(L)_{-j},\calC(Z,\Omega_0^\vee))^\bullet$
in the derived category $D_{\rm coh}^+(Z)$ is represented by the complex
$\calH om(\Lambda_Z(L)_{-j},\calC(Z,\Omega_0^\vee))^\bullet$
by \cite[Lem.3.6]{Ishida2}.
Hence we get the corollary by Theorem~\ref{Thm 5.1}.

The following corollary follows from Corollaries~\ref{Cor 2.12} and
\ref{Cor 5.2}.

%
%

\begin{Cor}
\label{Cor 5.3}
For each integer $0\leq j\leq r$, the Grothendieck dual of the
$d$-complex $\Omega_j(\p; Z)$ is quasi-isomorphic to
$\Omega_{r-j}(-\p; Z)[r]$.
\end{Cor}

We assume that $\Delta$ is a complete fan.
Then $Z$ is a complete variety.

Let $F^\bullet$ be a finite complex of coherent $\calO_Z$-modules.
By the Grothendieck duality theorem \cite[VI,Thm.3.3]{RD} applied for the
proper morphism $f_Z\cl Z\ra\Spec\Q$, we have a natural isomorphism
\begin{equation}
\H^i(Z, R\calH om(F,\calC(Z,\Omega_0^\vee))^\bullet)
\simeq\Hom(\H^{-i}(Z,F^\bullet),\Q)\;.
\end{equation}

%
%

\begin{Thm}
\label{Thm 5.4}
Assume that $\Delta$ is a complete fan.
Then the equality
\begin{equation}
\dim_\Q\H^i(Z,\Omega_j(\p; Z)) =\dim_\Q\H^{r-i}(Z,\Omega_{r-j}(-\p; Z))
\end{equation}
holds for any integers $i, j$.
\end{Thm}

\Proof
By Corollary~\ref{Cor 5.2}, this is a consequence of the
Grothendieck-Serre duality theorem.

We can also prove the equality directly as follows.
We have equalities
\begin{equation}
\dim_\Q\H^i(Z,\Omega_j(\p; Z)^\bullet) =
\dim_\Q\H^i(\Gamma(\ic_\p(\Delta)^\bullet)_{-j}
\end{equation}
and
\begin{equation}
\dim_\Q\H^{r-i}(Z,\Omega_{r-j}(\p; Z)^\bullet) =
\dim_\Q\H^{r-i}(\Gamma(\ic_{-\p}(\Delta)^\bullet)_{-r+j}
\end{equation}
by Lemma~\ref{Lem 3.6}.
Since $\ic_\p(\Delta)^\bullet$ is quasi-isomorphic to
$\D(\ic_{-\p}(\Delta))^\bullet$ by Corollary~\ref{Cor 2.12},
we get the equality by Proposition~\ref{Prop 2.5}.
\QED


\begin{thebibliography}{XXXX}


\bibitem[BBD]{BBD}
A.A.\ Beilinson, J.\ Bernstein and P.\ Deligne,
Faisceaux pervers,
Analyse et Topologie sur les Espaces Singuliers (I),
Ast\'erisque {\bf 100}, Soc.\ Math.\ France, 1982.

\bibitem[D1]{Danilov}
V.I.\ Danilov, The geometry of toric varieties,
Russian Math.\ Surveys {\bf 33}, (1978), 97-154.

\bibitem[D2]{Deligne}
P.\ Deligne, Cohomologie a supports propres,
Expose XVII, SGA4, Tome 3, Lecture Notes in Math.\ {\bf 305},
Springer-Verlag, Berlin, Heidelberg, New York, 1973.

\bibitem[DL]{DenefLoeser}
J.\ Denef and F.\ Loeser, Weights of exponential sums,
intersection cohomology, and Newton polyhedra, preprint.

\bibitem[EGA4]{EGA4}
A.\ Grothendieck and J.\ Dieudonn\'e,
\'El\'ements de G\'eom\'etrie Alg\'ebrique IV, Inst.\ Hautes
\'Etudes Sci.\ Publ.\ Math.\ {\bf 20}, {\bf 24}, {\bf 28}, {\bf 32},
(1964-1967).


\bibitem[GAGA]{GAGA}
J.-P.\ Serre,
G\'eom\'etrie Alg\'ebrique et G\'eom\'etrie Analytique,
Ann.\ Inst.\ Fourier {\bf 6}, (1956), 1-42.


\bibitem[GM1]{GM1}
M.\ Goresky and R.\ MacPherson,
Intersection homology theory,
Topology {\bf 19}, (1980), 135-162.

\bibitem[GM2]{GM2}
M.\ Goresky and R.\ MacPherson,
Intersection homology II,
Invent.\ math.\ {\bf 72}, (1983), 77-129.

\bibitem[G]{Greub}
W.H.\ Grueb,
Multilinear Algebra,
Die Grundlehren der mathematischen Wissenschaften in
Einzeldarstellungen, {\bf 136}, Springer-Verlag, Berlin,
Heidelberg, New York, 1967.

\bibitem[H]{Hartshorne}
R.\ Hartshorne,
On the de Rham cohomology of algebraic varieties,
Inst.\ Hautes \'Etudes Sci.\ Publ.\ Math.\ {\bf 45}, 1975.

\bibitem[I1]{Ishida1}
M.-N.\ Ishida,
Torus embeddings and dualizing complexes, Tohoku Math.\ J.\
{\bf 32}, (1980), 111-146.

\bibitem[I2]{Ishida2}
M.-N.\ Ishida,
Torus embeddings and de Rham complexes, in Commutative
Algebra and Combinatorics (M.\ Nagata and H.\ Matsumura, ed.),
Adv.\ Studies in Pure Math.\ {\bf 11}, Kinokuniya, Tokyo and
North-Holland, Amsterdam, New York, Oxford, 1987, 111-145.

\bibitem[I3]{Ishida3}
M.-N.\ Ishida,
Torus embeddings and algebraic intersection complexes, II,
preprint.

\bibitem[MV]{MacVilonen}
R.\ MacPherson and K.\ Vilonen,
Elementary construction of perverse sheaves,
Invent.\ math.\ {\bf 84}, (1986), 403-435.

\bibitem[O1]{Oda1}
T.\ Oda,
Convex Bodies and Algebraic Geometry,
An Introduction to the Theory of Toric Varieties,
Ergebnisse der Math.\ (3), {\bf 15}, Springer-Verlag,
Berlin, Heidelberg, New York, Tokyo, 1988.

\bibitem[O2]{Oda2}
T.\ Oda,
The algebraic de Rham theorem for toric varieties,
Tohoku Math.\ J.\ {\bf 45}, (1993), 231-247.

\bibitem[O3]{Oda3}
T.\ Oda,
Simple convex polytopes and the strong Lefschetz theorem,
J.\ of Pure and Appl.\ Alg.\ {\bf 71}, (1991), 265-286.

\bibitem[RD]{RD}
R.\ Hartshorne,
Residues and Duality,
Lecture Notes in Math.\ {\bf 20},
Springer-Verlag, Berlin, Heidelberg, New York, 1966.

\bibitem[SC]{AMRT}
A.\ Ash, D.\ Mumford, M.\ Rapoport and Y.\ Tai,
Smooth compactification of locally symmetric varieties,
Lie Groups: History, Frontiers and Applications IV,
Math.\ Sci.\ Press, Brookline, Mass., 1975.

\bibitem[S]{Stanley}
R.\ Stanley, Generalized $h$-vectors, intesection cohomology
of toric varieties, and related results, in Commutative
Algebra and Combinatorics (M.\ Nagata and H.\ Matsumura, ed.), Adv.\ Studies
in Pure Math.\ {\bf 11}, Kinokuniya, Tokyo and North-Holland,
Amsterdam, New York, Oxford, 1987, 187-213.

\bibitem[V]{Verdier}
J.-L.\ Verdier,
Dualit\'e dan la cohomologie des espaces localement compacts,
S\'eminar Bourbaki, {\bf 300}, (1965).


\end{thebibliography}
\end{document}